\pdfoutput=1
\RequirePackage{ifpdf}
\ifpdf % We~are running pdfTeX in pdf mode
\documentclass[pdftex]{sigma}
\else
\documentclass{sigma}
\fi

\usepackage{tikz}
\usetikzlibrary{arrows}

\usepackage{array}
\newcolumntype{r}{D{.}{.}{-1}}

\numberwithin{equation}{section}

\newtheorem{Theorem}{Theorem}[section]
\newtheorem*{Theorem*}{Theorem}
\newtheorem{Corollary}[Theorem]{Corollary}
\newtheorem{Lemma}[Theorem]{Lemma}
\newtheorem{Proposition}[Theorem]{Proposition}
 { \theoremstyle{definition}
\newtheorem{Definition}[Theorem]{Definition}

\newtheorem{Example}[Theorem]{Example}
\newtheorem{Remark}[Theorem]{Remark} }

\DeclareMathOperator{\trop}{trop}
\DeclareMathOperator{\ev}{ev}
\DeclareMathOperator{\val}{val}
\DeclareMathOperator{\vir}{vir}
\DeclareMathOperator{\mult}{mult}
\DeclareMathOperator{\Aut}{Aut}
\DeclareMathOperator{\Coef}{coef}
\DeclareMathOperator{\lo}{loop}
\DeclareMathOperator{\ft}{ft}

\newcommand {\PP}{{\mathbb P}}
\newcommand {\RR}{{\mathbb R}}
\newcommand {\ZZ}{{\mathbb Z}}

\newcommand {\NN}{{\mathbb N}}
\newcommand {\TT}{{\mathbb T}}
\newcommand {\CC}{{\mathbb C}}

\begin{document}
\allowdisplaybreaks

\newcommand{\arXivNumber}{1809.10659}

\renewcommand{\thefootnote}{}

\renewcommand{\PaperNumber}{046}

\FirstPageHeading

\ShortArticleName{Tropical Mirror Symmetry in Dimension One}

\ArticleName{Tropical Mirror Symmetry in Dimension One\footnote{This paper is a~contribution to the Special Issue on Enumerative and Gauge-Theoretic Invariants in honor of Lothar G\"ottsche on the occasion of his 60th birthday. The~full collection is available at \href{https://www.emis.de/journals/SIGMA/Gottsche.html}{https://www.emis.de/journals/SIGMA/Gottsche.html}}}

\Author{Janko B\"OHM~$^{\rm a}$, Christoph GOLDNER~$^{\rm b}$ and Hannah MARKWIG~$^{\rm b}$}
\AuthorNameForHeading{J.~B\"ohm, C.~Goldner and H.~Markwig}

\Address{$^{\rm a)}$~Fachbereich Mathematik,
TU Kaiserslautern, Postfach 3049, 67653 Kaiserslautern, Germany}
\EmailD{\href{mailto:boehm@mathematik.uni-kl.de}{boehm@mathematik.uni-kl.de}}

\Address{$^{\rm b)}$~Universit\"at T\"ubingen, Fachbereich Mathematik, 72076 T\"ubingen, Germany}
\EmailD{\href{mailto:christoph.goldner@math.uni-tuebingen.de}{christoph.goldner@math.uni-tuebingen.de}, \href{mailto:hannah@math.uni-tuebingen.de}{hannah@math.uni-tuebingen.de}}

\ArticleDates{Received January 24, 2022, in final form June 17, 2022; Published online June 25, 2022}

\Abstract{We prove a tropical mirror symmetry theorem for descendant Gromov--Witten invariants of the elliptic curve, generalizing the tropical mirror symmetry theorem for Hurwitz numbers of the elliptic curve, Theorem~2.20 in~[B\"ohm J., Bringmann K., Buchholz A., Markwig H., \textit{J.~Reine Angew. Math.} \textbf{732} (2017), 211--246, arXiv:1309.5893]. For the case of the elliptic curve, the tropical version of mirror symmetry holds on a fine level and easily implies the equality of the generating series of descendant Gromov--Witten invariants of the elliptic curve to Feynman integrals. To prove tropical mirror symmetry for elliptic curves, we investigate the bijection between graph covers and sets of monomials contributing to a~coefficient in a~Feynman integral. We also soup up the traditional approach in mathema\-tical physics to mirror symmetry for the elliptic curve, involving operators on a Fock space, to give a proof of tropical mirror symmetry for Hurwitz numbers of the elliptic curve. In this way, we shed light on the intimate relation between the operator approach on a bosonic Fock space and the tropical approach.}


\Classification{14J33; 14N35; 14T05; 81T18; 11F11; 14H30; 14N10; 14H52; 14H81}

\renewcommand{\thefootnote}{\arabic{footnote}}
\setcounter{footnote}{0}

\section{Introduction}

\subsection{Context: Tropical mirror symmetry of elliptic curves}

Mirror symmetry is a duality relation involving algebraic resp.\ symplectic varieties and their invariants. Its main motivation comes from string theory, but it is also at the base of many interesting developments in mathematics.
We focus on statements relating generating series of Gromov--Witten invariants of a variety $X$ with certain integrals on its mirror $X^\vee$.

Tropical geometry becomes a new tool to prove such relations, largely due to the well-known Gross--Siebert program, which aims at constructing new mirror pairs and providing an algebraic framework for SYZ-mirror symmetry \cite{GS06, GS07, SYZ}.

The philosophy how tropical geometry can be exploited is illustrated in the following triangle:

\begin{center}

\scalebox{0.95}{
%\centering
\begin{tikzpicture}[<->,>=stealth',shorten >=1pt,auto]
\coordinate (a) at (0,0);
\coordinate (b) at (-3.4,3.5);
\coordinate (c) at (3.4,3.5);
\node(1) at (a) {\begin{tabular}{c}tropical\\ GW-invariants\end{tabular}};
\node(2) at (b) {\begin{tabular}{c}Gromov--Witten\\ invariants\end{tabular}};
\node(3) at (c) {\begin{tabular}{c}Feynman\\ integrals\end{tabular}};
\path[every node/.style={}]
(1) edge node[left,sloped, anchor=center] {\begin{tabular}{c}Correspondence\\ Theorem\end{tabular}} (2)
(2) edge [dashed] node {Mirror symmetry} (3)
(3) edge node[right,sloped, anchor=center] {} (1);
%(3) edge node[right,sloped, anchor=center] {\begin{tabular}{c}Tropical\\ mirror symmetry\end{tabular}} (1);
\end{tikzpicture}
}

\end{center}

In many situations, correspondence theorems relating Gromov--Witten invariants resp.\ enumerative invariants to their tropical counterparts are known \cite{BBM10, CJM10, Mi03, NS06}. If we can relate the generating function of tropical invariants to integrals, we obtain a proof of the desired mirror symmetry relation using a detour via tropical geometry \cite{Gro09, Ove15}.

In \cite{BBBM13}, we investigated the triangle above for the case of Hurwitz numbers of the elliptic curve and Feynman integrals.
Correspondence theorems for Hurwitz numbers existed already, tropical Hurwitz numbers essentially count certain decorated graphs.
The mirror symmetry relation in this case was known, there is a proof in mathematical physics involving operators on a Fock space.

The tropical approach revealed that the relation holds on an even finer level: tropically, we can relate Feynman integrals and generating series of (labeled) tropical covers graph by graph and order by order. As a consequence, one obtains interesting new quasimodularity statements for graph generating series \cite{GM16}.

The mirror symmetry theorem (the top arrow) is an easy corollary of the more general tropical version. The tropical mirror symmetry theorem for Hurwitz numbers of the elliptic curve can be viewed as a support for the strategy of the Gross--Siebert program, or more generally for the philosophy of using tropical geometry as a tool in mirror symmetry.

\subsection[Part I: Generating series of tropical descendant Gromov--Witten invariants of E]
{Part I: Generating series of tropical descendant \\Gromov--Witten invariants of $\boldsymbol E$}

From Okounkov--Pandharipande's Gromov--Witten/Hurwitz (GW/H) correspondence \cite[Theorem~1]{OP06}, it is known that Hurwitz numbers are a special case of descendant Gromov--Witten invariants.

The central result of this article is a general tropical mirror symmetry theorem for elliptic curves, involving descendant Gromov--Witten invariants (Theorem \ref{thm-refinedmirrortrop}):

\begin{Theorem}
Generating series of tropical descendant Gromov--Witten invariants of an elliptic curve can be expressed in terms of Feynman integrals. In particular, they are quasimodular forms, also when restricted to a certain combinatorial type of source.
\end{Theorem}

Together with the suitable correspondence Theorem \ref{thm-corres} (see~\cite[Theorem~3.1.2]{CJMR16}) relating tropical descendant Gromov--Witten invariants to their counterparts in algebraic geometry, it can be applied to prove the mirror symmetry theorem for elliptic curves involving descendants (see Theorem \ref{thm-mirror},  \cite[resp.\ Theorem 6.1 and Proposition 6.7]{Lithesis} or \cite[Theorem 1.2(2) and Proposition 3.4]{Li11}).

The most important new tool of \cite{BBBM13} in the study of mirror symmetry for elliptic curves was a bijection between tropical covers (i.e., the above mentioned decorated graphs) satisfying fixed discrete data and sets of monomials contributing to a coefficient in a Feynman integral \cite[Theorem 2.30]{BBBM13}. For our purpose, we generalize this in two directions, both involving the source curves of the tropical covers in question:
\begin{enumerate}\itemsep=0pt
\item[$(a)$] we need to allow vertices of valency different from $3$, and
\item[$(b)$] we need to allow genus at vertices.
\end{enumerate}
The task in $(a)$ is a major extension of the bijection, formulated in Theorem \ref{thm-bij}, and will have further applications.
The task in $(b)$ involves the multiplicities with which covers are counted and is more a question of bookkeeping.

The tasks $(a)$ and $(b)$ are necessary, as psi-classes in the tropical world impose higher valency resp.\ higher genus on vertices, see, e.g., \cite{CGM20, KM06, MR20, MR08, Mi07}.
For example, for the tropical count of the descendant Gromov--Witten invariant
\begin{gather*}
\langle \tau_{2}(pt) \tau_{0}(pt) \tau_{0}(pt) \rangle_{2,3}^{E,3,\trop}
\end{gather*}
 (see Example \ref{ex-stationarydescendant}), the two graphs depicted in Figure \ref{fig-intro} are needed as combinatorial types of source curves for tropical covers.

\begin{figure}[h!]\centering

\tikzset{every picture/.style={line width=0.75pt}} %set default line width to 0.75pt

\begin{tikzpicture}[x=0.75pt,y=0.75pt,yscale=-1,xscale=1]
%uncomment if require: \path (0,784); %set diagram left start at 0, and has height of 784

%Shape: Ellipse [id:dp013580038362878022]
\draw (180.33,170.47) .. controls (180.33,159.42) and (196,150.47) .. (215.33,150.47) .. controls (234.66,150.47) and (250.33,159.42) .. (250.33,170.47) .. controls (250.33,181.51) and (234.66,190.47) .. (215.33,190.47) .. controls (196,190.47) and (180.33,181.51) .. (180.33,170.47) -- cycle ;
%Shape: Ellipse [id:dp7701962010957114]
\draw (53.33,173.8) .. controls (53.33,162.75) and (69,153.8) .. (88.33,153.8) .. controls (107.66,153.8) and (123.33,162.75) .. (123.33,173.8) .. controls (123.33,184.85) and (107.66,193.8) .. (88.33,193.8) .. controls (69,193.8) and (53.33,184.85) .. (53.33,173.8) -- cycle ;
%Shape: Ellipse [id:dp2042130485589423]
\draw (250.33,170.47) .. controls (250.33,159.42) and (266,150.47) .. (285.33,150.47) .. controls (304.66,150.47) and (320.33,159.42) .. (320.33,170.47) .. controls (320.33,181.51) and (304.66,190.47) .. (285.33,190.47) .. controls (266,190.47) and (250.33,181.51) .. (250.33,170.47) -- cycle ;
%Shape: Circle [id:dp1490748918822289]
\draw [fill={rgb, 255:red, 0; green, 0; blue, 0 } ,fill opacity=1 ] (248.58,170.47) .. controls (248.58,169.5) and (249.37,168.72) .. (250.33,168.72) .. controls (251.3,168.72) and (252.08,169.5) .. (252.08,170.47) .. controls (252.08,171.43) and (251.3,172.22) .. (250.33,172.22) .. controls (249.37,172.22) and (248.58,171.43) .. (248.58,170.47) -- cycle ;
%Shape: Circle [id:dp35774114475855245]
\draw [fill={rgb, 255:red, 0; green, 0; blue, 0 } ,fill opacity=1 ] (51.58,172.05) .. controls (51.58,171.08) and (52.37,170.3) .. (53.33,170.3) .. controls (54.3,170.3) and (55.08,171.08) .. (55.08,172.05) .. controls (55.08,173.02) and (54.3,173.8) .. (53.33,173.8) .. controls (52.37,173.8) and (51.58,173.02) .. (51.58,172.05) -- cycle ;
%Shape: Circle [id:dp6750115315428419]
\draw [fill={rgb, 255:red, 0; green, 0; blue, 0 } ,fill opacity=1 ] (71.58,192.05) .. controls (71.58,191.08) and (72.37,190.3) .. (73.33,190.3) .. controls (74.3,190.3) and (75.08,191.08) .. (75.08,192.05) .. controls (75.08,193.02) and (74.3,193.8) .. (73.33,193.8) .. controls (72.37,193.8) and (71.58,193.02) .. (71.58,192.05) -- cycle ;
%Shape: Circle [id:dp28707115245174253]
\draw [fill={rgb, 255:red, 0; green, 0; blue, 0 } ,fill opacity=1 ] (80.92,154.38) .. controls (80.92,153.42) and (81.7,152.63) .. (82.67,152.63) .. controls (83.63,152.63) and (84.42,153.42) .. (84.42,154.38) .. controls (84.42,155.35) and (83.63,156.13) .. (82.67,156.13) .. controls (81.7,156.13) and (80.92,155.35) .. (80.92,154.38) -- cycle ;
%Shape: Circle [id:dp6625924016292067]
\draw [fill={rgb, 255:red, 0; green, 0; blue, 0 } ,fill opacity=1 ] (281.83,150.47) .. controls (281.83,149.5) and (282.62,148.72) .. (283.58,148.72) .. controls (284.55,148.72) and (285.33,149.5) .. (285.33,150.47) .. controls (285.33,151.43) and (284.55,152.22) .. (283.58,152.22) .. controls (282.62,152.22) and (281.83,151.43) .. (281.83,150.47) -- cycle ;
%Shape: Circle [id:dp2998301737068383]
\draw [fill={rgb, 255:red, 0; green, 0; blue, 0 } ,fill opacity=1 ] (270.5,188.8) .. controls (270.5,187.83) and (271.28,187.05) .. (272.25,187.05) .. controls (273.22,187.05) and (274,187.83) .. (274,188.8) .. controls (274,189.77) and (273.22,190.55) .. (272.25,190.55) .. controls (271.28,190.55) and (270.5,189.77) .. (270.5,188.8) -- cycle ;

% Text Node
\draw (15.07,163.47) node [anchor=north west][inner sep=0.75pt] [font=\small] [align=left] {$g=1$};

\end{tikzpicture}

\caption{Two graphs that appear as combinatorial types of source curves of tropical covers contributing to $\langle \tau_{2}(pt) \tau_{0}(pt) \tau_{0}(pt) \rangle_{2,3}^{E,3,\trop}$. The left has a vertex of genus one, the right has a $4$-valent vertex.}\label{fig-intro}

\end{figure}

\subsection{Part II: Relation to the Fock space approach}

The traditional approach to mirror symmetry of an elliptic curve involves operators on Fock spaces.
There are two Fock spaces, a fermionic and a bosonic Fock space, and an isomorphism between them called the boson--fermion correspondence. The latter is usually viewed as the essence of mirror symmetry.
The generating function of Gromov--Witten invariants can be interpreted on the fermionic side, via the correspondence we then obtain an expression in terms of matrix elements on the bosonic Fock space, and the latter can be related to Feynman integrals \cite{KR, Lithesis, Li11, OP06}.

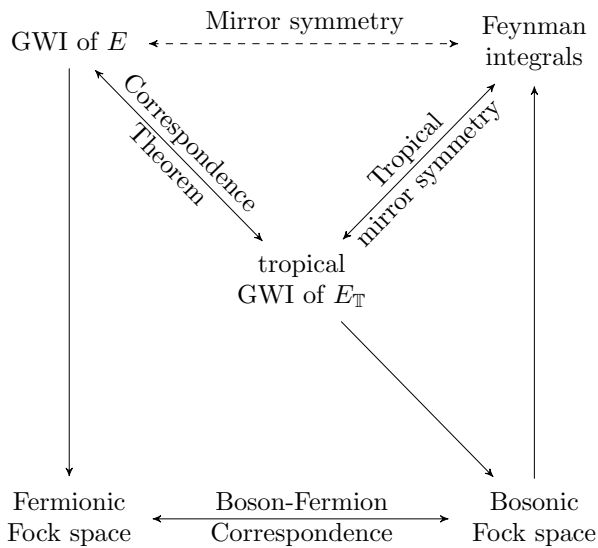
\begin{figure}[t]\centering
\scalebox{0.93}{
\begin{tikzpicture}[]
\coordinate (a) at (0,0);
\coordinate (b) at (-3.4,3.5);
\coordinate (c) at (3.4,3.5);
\coordinate (d) at (-3.4,-3.5);
\coordinate (e) at (3.4,-3.5);
\node(1) at (a) {\begin{tabular}{c}tropical\\ GWI of $E_{\TT}$\end{tabular}};
\node(2) at (b) {\begin{tabular}{c}GWI of $E$\end{tabular}};
\node(3) at (c) {\begin{tabular}{c}Feynman\\ integrals\end{tabular}};
\node(4) at (d) {\begin{tabular}{c}Fermionic \\ Fock space\end{tabular}};
\node(5) at (e) {\begin{tabular}{c}Bosonic \\ Fock space\end{tabular}};
\path[every node/.style={},<->,>=stealth',shorten >=1pt,auto]
(1) edge node[left,sloped, anchor=center] {\begin{tabular}{c}Correspondence\\ Theorem\end{tabular}} (2)
(2) edge [dashed] node {Mirror symmetry} (3)
(3) edge node[right,sloped, anchor=center] {\begin{tabular}{c}Tropical\\ mirror symmetry\end{tabular}} (1)
(4) edge node[left,sloped, anchor=center] {\begin{tabular}{c}boson--fermion\\ Correspondence\end{tabular}} (5);
\path[every node/.style={},->,>=stealth',shorten >=1pt,auto]
(5) edge node[left,sloped, anchor=center] {} (3)
(1) edge node[left,sloped, anchor=center] {} (5)
(2) edge node[left,sloped, anchor=center] {} (4);
\end{tikzpicture}
}

\caption{Boson--fermion correspondence and tropical geometry as a shortcut.}\label{fig-Overview}
\end{figure}

We reveal the close connection between the Fock space approach to mirror symmetry of elliptic curves and the new tropical approach.

Here, tropical geometry hands us a shortcut: by passing to tropical Gromov--Witten invariants on the tropical elliptic curve $E_\mathbb{T}$, we can directly relate the generating series of Gromov--Witten invariants of the elliptic curve $E$ to a matrix element on the bosonic Fock space (see Figure \ref{fig-Overview}), supporting the slogan ``tropicalization is bosonification'' from \cite{CJMR16}, and the intuition underlying the Gross--Siebert program that tropical geometry is a natural language in the context of mirror symmetry.

We soup up the traditional Fock space approach to give an alternative proof of the tropical mirror symmetry relation Theorem \ref{thm-refinedmirrortrop}, for simplicity restricting to the case of Hurwitz numbers.
The main ingredient is a version of Wick's theorem which encodes matrix elements in a~bosonic Fock space as weighted sums of graphs, which can then directly be related to tropical Hurwitz covers (see Theorem \ref{thm-fock}).

Since the first version of this paper appeared as preprint, other researchers have continued working on related topics. We would like to point out in particular a series of papers by Blomme who studies the enumerative geometry of line bundles over elliptic curves and generalized further to the enumerative geometry of abelian surfaces \cite{Blom21, Blom22}, and the papers on enumerative geometry of elliptic fibrations by Oberdieck and Pixton \cite{OP18, OP19}.

Section \ref{chap-mirror} focuses on tropical mirror symmetry and its direct proof via a bijection involving graph covers and monomials in a Feynman integral.
Again, the fact that tropical mirror symmetry holds on a fine level has interesting implications for graph summands of the generating series of descendant Gromov--Witten invariants. The quasimodularity of these graph summands is shown in Section \ref{sec-quasimod} relying on \cite[Theorem 6.1]{GM16}.
Section \ref{chap-fock} is devoted to the Fock space approach.

\section{Tropical mirror symmetry for elliptic curves}\label{chap-mirror}
\subsection{Descendant Gromov--Witten invariants}
Gromov--Witten invariants are virtually enumerative intersection numbers on moduli spaces of stable maps. Let $E$ be an elliptic curve. Gromov--Witten invariants of $E$ do not depend on its complex structure.
A {\it stable map} of degree $d$ from a curve of genus $g$ to $E$ with $n$ markings is a map $f\colon  C \to E$, where $C$ is a connected projective curve with at worst nodal singularities, and with $n$ distinct nonsingular marked points $x_1,\dots,x_n\in C$, such that $f_\ast([C]) = d[E]$ and~$f$ has a finite group of automorphism.
The moduli space of stable maps, denoted $\overline{\mathcal{M}}_{g,n} (E,d)$, is a proper Deligne--Mumford stack of virtual dimension $2g-2+n$ \cite{Beh97,BF97}.
The $i$th evaluation morphism is the map $\ev_i\colon  \overline{\mathcal{M}}_{g,n} (E,d) \to E$ that sends a point $[C, x_1, \dots, x_n, f]$ to $f(x_i)
\in E$.
The $i$th cotangent line bundle $\mathbb{L}_i \to \overline{\mathcal{M}}_{g,n} (E,d)$ is defined by a canonical identification of its fiber over a moduli point $(C, x_1, \dots, x_n, f)$
with the cotangent space $T^\ast_{x_i}(C)$. The first Chern class of the cotangent line bundle is called a {\it psi class} ($\psi_i= c_1(\mathbb{L}_i)$).

 \begin{Definition}
 Fix $g$, $n$, $d$ and let $k_1, \dots, k_n$ be non-negative integers with
 \begin{gather*}
 k_1+\dots+k_n = 2g-2.
  \end{gather*}
 The {\it stationary descendant Gromov--Witten invariant} $\langle \tau_{k_1}(pt) \cdots \tau_{k_n}(pt) \rangle_{g,n}^{E,d}$ is defined by
 \begin{equation*}
 \langle \tau_{k_1}(pt) \cdots \tau_{k_n}(pt) \rangle_{g,n}^{E,d} = \int_{[\overline{\mathcal{M}}_{g,n} (E,d)]^{\rm vir}} \prod_{i=1}^n \ev_i^\ast(pt) \psi_i^{k_i},
 \end{equation*}
 where $pt$ denotes the class of a point in $E$.
 \end{Definition}

In Section \ref{chap-fock}, we use degeneration techniques to relate the proof of mirror symmetry for elliptic curves in mathematical physics to the tropical approach. For this purpose, we also need to introduce {\it relative Gromov--Witten invariants}:
they are constructed using moduli spaces of {\it relative stable maps} $\overline{\mathcal{M}}_{g,n} \big(\PP^1, \mu,\nu,d\big)$, where part of the data specified are the ramification profiles $\mu$ and $\nu$ which we fix over $0$ resp.\ $\infty\in \PP^1$. The preimages of $0$ and $\infty$ are marked. A~detailed discussion of spaces of relative stable maps to $\PP^1$ and their boundary is not necessary for our purpose, we refer to \cite{Vak08}. We use operator notation and denote
\begin{equation*}
 \langle \mu| \tau_{k_1}(pt) \cdots \tau_{k_n}(pt)|\nu \rangle_{g, n}^{\PP^1,d} = \int_{[\overline{\mathcal{M}}_{g,n} (\PP^1, \mu,\nu,d)]^{\vir}} \prod_{i=1}^n \ev_i^\ast(pt) \psi_i^{k_i}.
\end{equation*}
One can allow source curves to be disconnected, and introduce {\it disconnected Gromov--Witten invariants}. We will add the superscript $\bullet$ anytime we wish to refer to the disconnected theory.

\begin{Remark}\label{rem-ki1covers}
It follows from the GW/H correspondence \cite[Theorem 1]{OP06} that a stationary descendant Gromov--Witten invariant with $k_i=1$ for all $i$ is a Hurwitz number counting covers of the resp. degree and genus and with $n$ fixed simple branch points.
\end{Remark}

\subsection{Tropical descendant Gromov--Witten invariants}

An \emph{abstract tropical} \emph{cur-ve} is a connected metric graph $\Gamma$, such that edges leading to leaves (called \emph{ends}) have infinite length, together with a genus function $g\colon \Gamma\rightarrow \ZZ_{\geq 0}$ with finite support. Locally around a point $p$, $\Gamma$ is homeomorphic to a star with $r$ halfrays.
The number $r$ is called the \emph{valence} of the point $p$ and denoted by $\val(p)$. We identify the vertex set of $\Gamma$ as the points where the genus function is nonzero, together with points of valence different from $2$. The vertices of valence greater than $1$ are called \textit{inner vertices}. Besides \emph{edges}, we introduce the notion of \emph{flags} of $\Gamma$. A flag is a pair $(V,e)$ of a vertex $V$ and an edge $e$ incident to it ($V\in \partial e$). Edges that are not ends are required to have finite length and are referred to as \emph{bounded} or \textit{internal} edges.

A \emph{marked tropical curve} is a tropical curve whose leaves are labeled. An isomorphism of a~tropical curve is an isometry respecting the leaf markings and the genus function. The \emph{genus} of a~tropical curve $\Gamma$ is given
\begin{gather*}
g(\Gamma) = h_1(\Gamma)+\sum_{p\in \Gamma} g(p).
\end{gather*}
A curve of genus $0$ is called \emph{rational} and a curve satisfying $g(v)=0$ for all $v$ is called \emph{explicit}. The \emph{combinatorial type} is the equivalence class of tropical curves obtained by identifying any two tropical curves which differ only by edge lengths.

A \emph{tropical cover} $\pi\colon \Gamma_1\rightarrow \Gamma_2$ is a surjective harmonic map of metric graphs in the sense of~\cite[Section 2]{ABBR1}. The map $\pi$ is piecewise integer affine linear, the slope of $\pi$ on a flag or edge $e$ is a~nonnegative integer called the \emph{expansion factor} $\omega(e)\in \NN$.

The expansion factor of $e$ can be $0$ only if $e$ is an end. We fix the convention that the ends marked $1,\dots,n$ are the ones with expansion factor $0$.

For a point $v\in \Gamma_1$, the \emph{local degree of $\pi$ at $v$} is defined as follows. Choose a flag $f'$ adjacent to $\pi(v)$, and add the expansion factors of all flags $f$ adjacent to $v$ that map to $f'$:
\begin{equation*}
d_v=\sum_{f\mapsto f'} \omega(f).
\end{equation*}
We define the \textit{harmonicity} or \textit{balancing condition} to be the fact that for each point $v\in \Gamma_1$, the local degree at $v$ is well defined (i.e., independent of the choice of $f'$).

The \emph{degree} of a tropical cover is the sum over all local degrees of preimages of a point $a$, $d=\sum_{p\mapsto a} d_p$ (here, we consider the map locally around a vertex of the source graph). By the balancing condition, this definition does not depend on the choice of $a\in \Gamma_2$. For a flag $f$ of the image graph $\Gamma_2$, let $\mu_f$ be the partition of expansion factors of the flags of the source graph $\Gamma_1$ mapping onto $f$. We call $\mu_f$ the \emph{ramification profile} above $f$.

The tropical projective line, $\PP^1_{\TT}$, equals $\RR\cup \{\pm \infty\}$, a (nondegenerate) tropical elliptic curve~$E_{\TT}$ is a circle with a fixed length.

\begin{Definition}[psi- and point conditions]\label{def-PsiAndPointConditions}
We say that a tropical cover $\pi\colon \Gamma_1\rightarrow \Gamma_2$ with a~marked end $i$ {\it satisfies a psi-condition} with power $k$ at $i$, if the vertex $V$ to which the marked end $i$ is adjacent has valency $k+3-2g(V)$. We say $\pi\colon \Gamma_1\rightarrow \Gamma_2$ {\it satisfies the point conditions} $p_1,\dots,p_n\in\Gamma_2$ if
\begin{gather*}
\lbrace\pi(1),\dots,\pi(n)\rbrace=\lbrace p_1,\dots,p_n \rbrace.
\end{gather*}
\end{Definition}

Fix $g$, $n$, $d$ and let $k_1, \dots, k_n$ be non-negative integers with
\begin{gather*}
k_1+\dots+k_n = 2g-2.
\end{gather*}
Let $\pi\colon \Gamma\rightarrow E_{\TT}$ be a tropical cover such that $\Gamma$ is of genus $g$ and has $n$ marked ends.
Fix~$n$ distinct points $p_1,\dots,p_n\in E_{\TT}$.
Assume that at the marked end $i$, a psi-condition with power~$k_i$ is satisfied, and that the point conditions are satisfied.
The marked ends must be adjacent to different vertices, since they satisfy different point conditions. It follows from an Euler characteristic argument incorporating the valencies imposed by the psi-conditions that $\Gamma$ has exactly $n$ vertices, each adjacent to one marked end.

Locally at the marked end $i$, the cover sends the vertex to an interval consisting of two flags~$f$ and $f'$. We define the {\it local vertex multiplicity} $\mult_i(\pi)$ to be a combinatorial factor times a~one-point relative descendant Gromov--Witten invariant:
\begin{equation} \mult_i(\pi)= \langle \mu_f| \tau_{k_i}(pt)|\mu_{f'} \rangle_{g_i, 1}^{\PP^1,d_i},\label{eq-localmult}\end{equation}
where $g_i$ denotes the genus of the vertex adjacent to the marked end $i$, $d_i$ its local degree, and~$\mu_f$ resp.\ $\mu_{f'}$ the ramification profiles above the two flags of the image interval.

We define the multiplicity of $\pi$ to be
\begin{equation}
\frac{1}{|\Aut(\pi)|}\cdot \prod_i \mult_i(\pi)\cdot \prod_e \omega(e).\label{eq-mult}
\end{equation}

Note that all ends of a tropical cover of $E_{\TT}$ are contracted ends, with image points the points~$p_i$ we fix as conditions in $E_{\TT}$.

\begin{Definition}[tropical stationary descendant Gromov--Witten invariant of $E_{\TT}$]
For $g$, $n$, $d$, $k_1,\dots,k_n$ as above, define the \emph{tropical stationary descendant Gromov--Witten invariant }
\begin{gather*}
\langle \tau_{k_1}(pt) \cdots \tau_{k_n}(pt) \rangle_{g,n}^{E,d,\trop}
\end{gather*}
to be the weighted count of tropical genus $g$ degree $d$ covers of $E_{\TT}$ with $n$ marked points satisfying point and psi-conditions as above, each counted with its multiplicity as defined in (\ref{eq-mult}).

\end{Definition}

Note that the metric structure of the source curves of covers contributing to a tropical des\-cendant Gromov--Witten invariant is implicit in the metric data of $E_{\TT}$ and the chosen point conditions. We can thus neglect length data in the source curve.

\begin{Example}\label{ex-stationarydescendant}
As an example, fix three different points $p_1$, $p_2$, $p_3$ on $E_{\TT}$ and let $d=3$, $g=2$, $k_1=2$, $k_2=0$, $k_3=0$. Note that $\sum_i k_i=2g-2$ is satisfied. We list all covers contributing to $\langle \tau_{2}(pt) \tau_{0}(pt) \tau_{0}(pt) \rangle_{2,3}^{E,3,\trop}$ in Figure \ref{fig-Example1} below. Figure \ref{fig-Example1} shows schematic representations of the source curves of all covers contributing, where we assume that the top vertex of each such representation is mapped to $p_1$, the right vertex is mapped to~$p_2$ and the left one is mapped to~$p_3$. This convention gives us one choice out of $3!$ choices of an order of labeled vertices of the source curve mapping to $p_1$, $p_2$, $p_3$ on $E_{\TT}$. A green number indicates that there is a nonzero genus $g_i$ at a vertex $i$. The other numbers are the weights of the edges that are greater than $1$. Note that the valency of a vertex $i$ is given by $k_i+3-2g_i$ when taking the contracted ends into account. When neglecting marked ends, the underlying graph is either a figure $8$ or a loop (see Example~\ref{ex-underlyingFeynmanGraph}). In each case, every loop is mapped to $E_\TT$. When we draw a curl in an edge, it means that the edge is mapped once around $E_\TT$.

\begin{figure}[h!]\centering
\def\svgwidth{395pt}
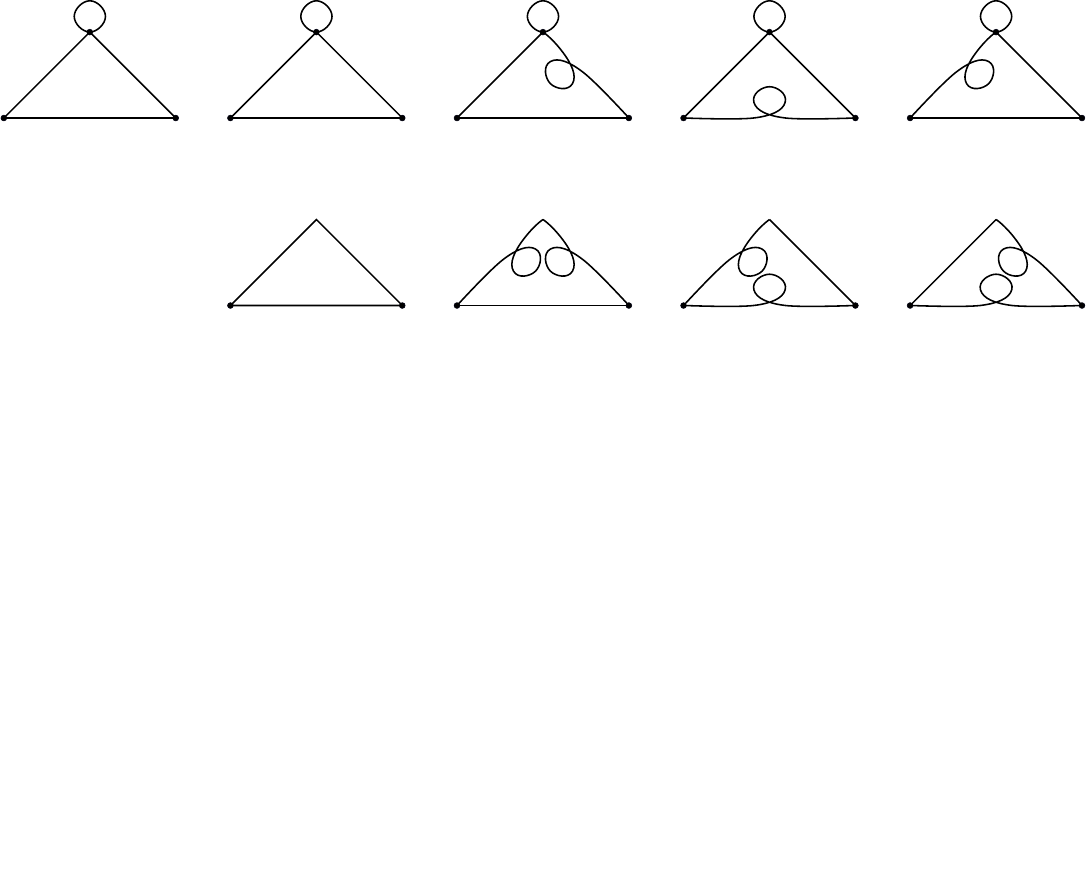
\caption{Schematic representations of source curves.}\label{fig-Example1}
\end{figure}

\noindent The multiplicity with which each curve is contributing is give by \eqref{eq-mult}. The local multiplicities $\mult_i(\pi)$, which are one-point relative descendant Gromov--Witten invariants, can be calculated explicitly using the one-point series \eqref{eq-onepointseries}. Each entry of the tabular below corresponds to one source curve of Figure \ref{fig-Example1} in the obvious way. An entry is the multiplicity of the corresponding cover $\pi$, where the first factor equals $|\Aut(\pi)|^{-1}$, the second factor equals $\prod_i \mult_i(\pi)$ and the third factor equals $\prod_e \omega(e)$.

{\setlength{\tabcolsep}{1.5em}
\renewcommand{\arraystretch}{1.3}

\begin{center}
\begin{tabular}{ c | c | c | c | c }
\hline
$1 \cdot 1 \cdot 2$&$1 \cdot 1 \cdot 8$&$1 \cdot 1 \cdot 1$&$1 \cdot 1 \cdot 1$&$1\cdot 1 \cdot 1$\\
\hline
$1 \cdot 1 \cdot 1$&$1 \cdot \frac{17}{24} \cdot 27$&$1 \cdot \frac{1}{24} \cdot 1$&$1 \cdot \frac{1}{24} \cdot 1$&$1\cdot \frac{1}{24} \cdot 1$\\
\hline
$1 \cdot \frac{1}{24} \cdot 1$&$1 \cdot \frac{1}{24} \cdot 1$&$1 \cdot \frac{1}{24} \cdot 1$&$1 \cdot 1 \cdot 4$&$1\cdot 1 \cdot 4$\\
\hline
$1 \cdot \frac{1}{24} \cdot 1$&$1 \cdot \frac{1}{24} \cdot 1$&$1 \cdot \frac{1}{24} \cdot 1$&$1 \cdot 1 \cdot 1$&$1\cdot 1 \cdot 1$\\
\hline
$1 \cdot 1 \cdot 1$&$1 \cdot 1 \cdot 1$&$1 \cdot 1 \cdot 1$&&\\
\hline
\end{tabular}
\end{center}}

\noindent Summing over all entries and considering the factor $3!$ yields
\begin{gather*}
\langle \tau_{2}(pt) \tau_{0}(pt) \tau_{0}(pt) \rangle_{2,3}^{E,3,\trop} =3!\cdot\frac{93}{2} =279.
\end{gather*}
\end{Example}

\begin{Theorem}[correspondence Theorem~I]\label{thm-corres}
A stationary descendant Gromov--Witten invariant of $E$ coincides with its tropical counterpart:
\begin{gather*}
\langle \tau_{k_1}(pt) \cdots \tau_{k_n}(pt) \rangle_{g,n}^{E,d}=\langle \tau_{k_1}(pt) \cdots \tau_{k_n}(pt) \rangle_{g,n}^{E,d,\trop}.
\end{gather*}
\end{Theorem}
For a proof, see \cite[Theorem~3.2.1]{CJMR16}.

To define tropical relative stationary descendant Gromov--Witten invariants of $\PP^1_{\TT}$, we fix two partitions $\mu$ and $\nu$ of the degree $d$. We consider tropical covers of $\PP^1_{\TT}$ such that the ramification profile over $-\infty$ equals $\mu$ and the ramification profile over $\infty$ equals $\nu$.
That is, in addition to the contracted ends that we use to impose point conditions, the source curve $\Gamma$ has $\ell(\mu)+\ell(\nu)$ marked ends which map to $\pm\infty$ with expansion factors imposed by $\mu$ and $\nu$. We assume that a~cover $\pi\colon \Gamma\rightarrow \PP^1_{\TT}$ meets point and psi-conditions as above. Local vertex multiplicities are defined as in equation~(\ref{eq-localmult}), and the multiplicity is
\begin{equation*}
\frac{1}{|\Aut(\pi)|}\cdot \prod_i \mult_i(\pi)\cdot \prod_e \omega(e),
%\label{eq-multrel}
\end{equation*}
where the last product goes over the bounded edges $e$ of $\Gamma$.
Tropical relative stationary descendant Gromov--Witten invariants of $\PP_{\TT}^1$, $\langle \mu| \tau_{k_1}(pt) \cdots \tau_{k_n}(pt)|\nu \rangle_{g, n}^{\PP^1,d, \trop}$, are defined as counts of tropical covers with the expansion factors of the unmarked ends imposed by $\mu$ and $\nu$ and satisfying the point and psi-conditions, counted with their multiplicity (see \cite[Definition~3.1.1]{CJMR16}).

\begin{Example}\label{ex-relativestationarydescendant}
Choose three different points $p_1$, $p_2$, $p_3$ on $E_{\TT}$ and let $d=3$, $g=2$, $k_1=2$, $k_2=0$, $k_3=0$ be as in Example \ref{ex-stationarydescendant}. Let $p_0$ be a base point on $E_{\TT}$ such that $p_0$, $p_1$, $p_2$, $p_3$ are ordered this way on $E_{\TT}$. Consider the source curve of a cover $\pi$ of $E_{\TT}$ depicted in the upper left corner of Figure \ref{fig-Example1} and cut it along $\pi^{-1}(p_0)$. Stretching the cut edges to infinity yields the cover shown below (we let $i$ be mapped to $p_i$). Note that this is a cover $\pi'$ to $\PP_{\TT}^1$ that contributes to $\langle (2,1)|\tau_{2}(pt) \tau_{0}(pt) \tau_{0}(pt)| (2,1) \rangle_{0,3}^{\PP^1,3,\trop}$.

\begin{figure}[h!]
\centering
\def\svgwidth{200pt}
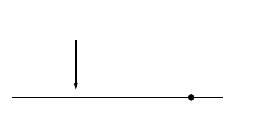
\end{figure}
\end{Example}

\begin{Theorem}[correspondence Theorem II]\label{thm-corresrel}
A relative $($stationary$)$ descendant Gromov--Witten invariant of $\PP^1$ coincides with its tropical counterpart:
\begin{gather*}
\langle \mu| \tau_{k_1}(pt) \cdots \tau_{k_n}(pt)|\nu \rangle_{g, n}^{\PP^1,d}=\langle \mu| \tau_{k_1}(pt) \cdots \tau_{k_n}(pt)|\nu \rangle_{g, n}^{\PP^1,d, \trop}.
\end{gather*}
\end{Theorem}
For a proof, see \cite[Theorem 3.1.2]{CJMR16}.

Just as before, we can also in the tropical world allow source curves to be disconnected, and add the superscript $\bullet$ in the notation.

\begin{Remark}[leaking]\label{rem-leaking}
We can tweak the definition of tropical covers of $E_{\TT}$ (resp. $\PP^1_{\TT}$) satisfying point and psi-conditions as follows: fix a direction for the target curve and specify for each end~$i$ of the source curve an integer $l_i$. Change the balancing condition in such a way that for the two flags $f_1$ and $f_2$ adjacent to $\pi(i)\in\lbrace p_1,\dots,p_n\rbrace$ (where we chose the notation to match the direction), the local degrees are not equal but differ by $l_i$:
\begin{gather*}
\sum_{f'\mapsto f_1} \omega(f')= \sum_{f''\mapsto f_2} \omega(f'')-l_i.
\end{gather*}
We call such covers \emph{leaky tropical covers}. Leaky tropical covers show up as floor diagrams representing counts of tropical curves in toric surfaces resp.\ in $\mathbb{P}^1$-bundles over $E$ (see, e.g., \cite{AB14, Blom21, Blom22, BM08, FM09}). We introduce them here, since they can be treated in terms of Feynman integrals analogously to their balanced versions.
\end{Remark}

\subsection{Feynman integrals}

\begin{Definition}[the propagator and the $\mathcal{S}$-function]\label{def-prop}
 We define the \emph{propagator} as a (formal) series in $x$ and $q$:
\begin{gather*}
 P(x,q)= \sum_{w=1}^\infty w\cdot x^{w} +\sum_{a=1}^\infty \bigg(\sum_{w|a}w \big(x^{w}+ x^{-w}\big)\bigg)q^{a}
\end{gather*}
and the $\mathcal{S}$-function as series in $z$:
\begin{gather*}
\mathcal{S}(z)=\frac{\sinh(z/2)}{z/2}.
\end{gather*}

We also consider another formal series in $q$ (which should be viewed as the propagator for loop edges):
\begin{gather*}
P^{\lo}(q)= \sum_{a=1}^\infty \biggl(\sum_{w|a} w\biggl) q^a.
\end{gather*}
\end{Definition}

\begin{Definition}[Feynman graphs]\label{def-Feynmangraphs}
Fix $n>1$. A \emph{Feynman graph} is a (non-metrized) graph~$\Gamma$ without ends with $n$ vertices which are labeled $x_1,\dots,x_n$ and with labeled edges $q_1,\dots,q_r$. By~convention, we assume that $q_1,\dots,q_s$ are loop edges and $q_{s+1},\dots,q_r$ are non-loop edges.
\end{Definition}
We do not fix the number of edges for a Feynman graph, the index $r$ can vary from graph to graph. We always use the letter $r$ for the number of edges in a fixed Feynman graph $\Gamma$.

\begin{Example}\label{ex-underlyingFeynmanGraph}
Recall Example \ref{ex-stationarydescendant}, where we provided all covers contributing to
\begin{gather*}
\langle \tau_{2}(pt) \tau_{0}(pt) \tau_{0}(pt) \rangle_{2,3}^{E,3,\trop}.
\end{gather*}
We can label their source curves, turning them into Feynman graphs, see Figure \ref{fig-Example12}.

\begin{Definition}[Feynman integrals] \label{def-Feynman} Let $\Gamma$ be a Feynman graph.
Let $\Omega$ be an order of the $n$ vertices of $\Gamma$.

For $k>s$, denote the vertices adjacent to the (non-loop) edge $q_k$ by $x_{k^1}$ and $x_{k^2}$, where we assume $x_{k^1}<x_{k^2}$ in $\Omega$.

For integers $l_1,\dots,l_n$, we define the \emph{Feynman integral} for $\Gamma$ and $\Omega$ to be
\begin{gather*}
I^{l_1,\dots,l_n}_{\Gamma,\Omega}(q)=\Coef_{[x_1^{l_1}\cdots x_n^{l_n}]} \prod_{k=1}^s P^{\lo}(q) \cdot \prod_{k=s+1}^{r} P\bigg(\frac{x_{k^1}}{x_{k^2}},q\bigg)
\end{gather*}
and the \emph{refined Feynman integral} to be\vspace{-.5ex}
\begin{gather*}
I^{l_1,\dots,l_n}_{\Gamma,\Omega}(q_1,\dots,q_r)=\Coef_{[x_1^{l_1}\cdots x_n^{l_n}]}\prod_{k=1}^s P^{\lo}(q_k) \prod_{k=s+1}^{r} P\bigg(\frac{x_{k^1}}{x_{k^2}},q_k\bigg).
\end{gather*}

Finally, we set\vspace{-.5ex}
\begin{gather*}
I^{l_1,\dots,l_n}_{\Gamma}(q)=\sum_\Omega I^{l_1,\dots,l_n}_{\Gamma,\Omega}(q),\vspace{-1ex}
\end{gather*}
where the sum goes over all $n!$ orders of the vertices of $\Gamma$, and\vspace{-.5ex}
\begin{gather*}
I^{l_1,\dots,l_n}_{\Gamma}(q_1,\dots,q_r)= \sum_\Omega I^{l_1,\dots,l_n}_{\Gamma,\Omega}(q_1,\dots,q_r).
\end{gather*}

If we drop the superscript $l_1,\dots,l_n$ in the notations above, then this stands for $l_i=0$ for all~$i$.
\end{Definition}\vspace{-1ex}
\begin{figure}[h!]
\centering
\def\svgwidth{395pt}
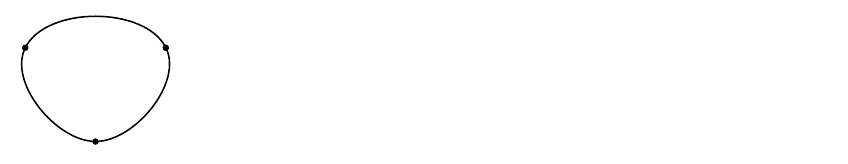\vspace{-3ex}
%\caption{??}
\caption{}
\label{fig-Example12}
\end{figure}
\end{Example}

If we assume $|x|<1$ to express the (in $q$) constant coefficient of the (non-loop) propagator (i.e., the first sum appearing in the propagator series in Definition \ref{def-prop}) as the rational function $\frac{x^2}{( x^2-1)^2}$ (using geometric series expansion), we can view the series from which we take the $x_1^{l_1}\cdots x_n^{l_n}$-coefficient in the Feynman integral above as a function on a Cartesian product of elliptic curves. If $l_i=0$ for all $i$, the Feynman integral then becomes a path integral in complex analysis (see~\cite[Definition~2.5 and equation~(2.4)]{BBBM13}).
Note that using the change of coordinates $x={\rm e}^{{\rm i}\pi u}$ the (non-loop) propagator has the following nice form
\begin{gather*}
P(u,q)=-\frac{1}{4\pi^2}\wp(u,q)-\frac{1}{12}E_2\big(q^2\big)
\end{gather*}
in terms of the Weierstra\ss{}-P-function $\wp$ and the Eisenstein series
\begin{gather*}
E_2(q):=1-24\sum_{d=1}^\infty \sigma(d)q^d.
\end{gather*}
Here, $\sigma$ denotes the sum-of-divisors function $\sigma(d)=\sum_{m|d}m$.
The variable $q$ above should be considered as a coordinate of the moduli space of elliptic curves, the variable $u$ as the complex coordinate of a fixed elliptic curve. (More precisely, $q={\rm e}^{2{\rm i}\pi \tau}$, where $\tau \in \CC$ is the parameter in the upper half plane in the well-known definition of the Weierstra\ss{}-P-function.)

\begin{Definition}[Feynman integrals with vertex contributions]\label{def-Feynmanvertex}
Let $\Gamma$ be a Feynman graph, and equip it with an additional genus function $\underline{g}$ associating a nonnegative integer $g_i$ to every vertex $x_i$. Let $\Omega$ be an order of the $n$ vertices of $\Gamma$.
We adapt our notion of propagators from Definitions~\ref{def-prop} and~\ref{def-Feynman} to include vertex contributions:
for non-loop edges, we set
\begin{gather*}
\begin{split}
&\tilde{P}\bigg(\frac{x_{k^1}}{x_{k^2}},q\bigg) =
\sum_{w=1}^\infty \mathcal{S}(w z_{k^1} )\mathcal{S}(w z_{k^2}) \cdot w \cdot \bigg(\frac{x_{k^1}}{x_{k^2}}\bigg)^w
\\
&\hphantom{\tilde{P}\bigg(\frac{x_{k^1}}{x_{k^2}},q\bigg) =} + \sum_{a=1}^\infty \bigg(\sum_{w|a}
\mathcal{S}(w z_{k^1} )\mathcal{S}(w z_{k^2})\cdot w \cdot \bigg(\bigg(\frac{x_{k^1}}{x_{k^2}}\bigg)^w+ \bigg(\frac{x_{k^2}}{x_{k^1}}\bigg)^w\bigg)\bigg)\cdot q^a.
\end{split}
\end{gather*}
For loop-edges connecting the vertex $x_{k^1}$ to itself, we set
\begin{gather*}
\tilde{P}^{\lo}(q)= \sum_{a=1}^\infty\bigg(\sum_{w|a}\mathcal{S}(w z_{k^1} )^2\cdot w\bigg) q^a.
\end{gather*}
The variables $z_{k_i}$ are new variables introduced for each vertex in order to take care of the genus contribution.

We define the \emph{Feynman integral with vertex contributions} for $\Gamma$, $\underline{g}$ and $\Omega$ to be
\begin{gather*}
I^{l_1,\dots,l_n}_{\Gamma,\underline{g},\Omega}(q)= \Coef_{[z_1^{2g_1}\cdots z_n^{2g_n}]} \Coef_{[x_1^{l_1}\cdots x_n^{l_n}]} \prod_{i=1}^n\frac{ 1}{\mathcal{S}(z_i)} \prod_{k=1}^s \tilde{P}^{\lo}(q) \prod_{k=s+1}^{r} \tilde{P}\bigg(\frac{x_{k^1}}{x_{k^2}},q\bigg)
\end{gather*}
and the \emph{refined Feynman integral with vertex contributions}
\begin{gather*}
I^{l_1,\dots,l_n}_{\Gamma,\underline{g},\Omega}(q_1,\dots,q_r)= \Coef_{[z_1^{2g_1}\cdots z_n^{2g_n}]} \Coef_{[x_1^{l_1}\cdots x_n^{l_n}]} \prod_{i=1}^n\frac{ 1}{\mathcal{S}(z_i)} \prod_{k=1}^s \tilde{P}^{\lo}(q_k) \!\!\prod_{k=s+1}^{r}\!\! \tilde{P}\bigg(\frac{x_{k^1}}{x_{k^2}},q_k\bigg) .
\end{gather*}

Again, we set
\begin{gather*}
I^{l_1,\dots,l_n}_{\Gamma,\underline{g}}(q)=\sum_\Omega I^{l_1,\dots,l_n}_{\Gamma,\underline{g},\Omega}(q),
\end{gather*}
where the sum goes over all $n!$ orders of the vertices, and
\begin{gather*}
I^{l_1,\dots,l_n}_{\Gamma,\underline{g}}(q_1,\dots,q_r)= \sum_\Omega I^{l_1,\dots,l_n}_{\Gamma,\underline{g},\Omega}(q_1,\dots,q_r).
\end{gather*}

Also here, dropping the superscript $l_1,\dots,l_n$ refers to the case $l_i=0$ for all $i$.
\end{Definition}

\subsection{(Tropical) mirror symmetry for elliptic curves}
\begin{Theorem}[mirror symmetry for $E$]\label{thm-mirror}
Fix $g\geq 2$, $n\geq 1$ and $k_1, \dots, k_n\geq 1$ satisfying $k_1+\dots+k_n = 2g-2$.

We can express the series of descendant Gromov--Witten invariants of $E$ in terms of Feynman integrals:
\begin{gather*}
\sum_{d\geq 1}\langle \tau_{k_1}(pt) \cdots \tau_{k_n}(pt) \rangle_{g,n}^{E,d}q^d = \sum_{(\ft(\Gamma),\underline{g})} \frac{1}{|\Aut(\ft(\Gamma),{\underline{g}})|}I_{\Gamma,\underline{g}}(q),
\end{gather*}
 where $\Gamma$ is a Feynman graph $($see Definition $\ref{def-Feynmangraphs})$ with a genus function $\underline{g}$, such that the vertex~$x_i$ has genus $g_i$ and valency $k_i+2-2g_i$, and such that $h^1(\Gamma)+\sum g_i=g$, and where we consider automorphisms of unlabeled graphs $(\ft$ is the forgetful map that forgets all labels of a Feynman graph $\Gamma$, see Definition~$\ref{def-LabeledTropCover})$ that are required to respect the genus function.
\end{Theorem}

A version of Theorem \ref{thm-mirror} is proved in \cite[Proposition 6.7]{Lithesis} (resp.\ \cite[Proposition 3.4]{Li11}) using the Fock space approach common in mathematical physics to which we relate the tropical approach in Section \ref{chap-fock}. In our approach, Theorem \ref{thm-mirror} becomes an easy corollary obtained by combining our Tropical mirror symmetry Theorem \ref{thm-refinedmirrortrop} with the correspondence Theorem \ref{thm-corres}.

\begin{Example}[automorphisms]\label{ex-automorphisms}
Consider the middle Feynman graph of Example~\ref{ex-underlyingFeynmanGraph}, denote it by $\Gamma$ and let its genus function be $\underline{g}=0$, i.e., there is no genus at the vertices. The automorphisms appearing in Theorem \ref{thm-mirror} are automorphisms respecting the underlying graph structure and the genus function of $(\Gamma ,\underline{g})$. In other words, we forget the labels of $\Gamma$ before determining its automorphisms. In case of $\Gamma$ as above, the automorphism group is $\ZZ_2\times\ZZ_2\times\ZZ_2$, because we can exchange the edges $q_1$ and $q_2$ (see Example \ref{ex-underlyingFeynmanGraph}) which gives a factor of $\ZZ_2$, we can exchange the edges $q_3$ and $q_4$ and we can exchange the vertices $x_2$ and $x_3$ in such a way that the edge $q_1$ maps to $q_3$ and the edge $q_2$ maps to $q_4$, see also the left side of Figure \ref{fig-Example13}.

In Section \ref{sec-quasimod}, we deal with unlabeled tropical covers, but with fixed \textit{order}. That is, we fix which end $i$ maps to which point $p_j$ on the elliptic curve. In such a case, on the Feynman integral side, we deal with automorphisms of the underlying Feynman graph with vertex labels (see Corollary \ref{cor-partiallyLabeled}). If we choose $(\Gamma,\underline{g})$ as above, then the automorphism group of the graph with vertex labels is $\ZZ_2\times\ZZ_2$ since we cannot exchange the vertices $x_2$ and $x_3$ anymore, they are now distinguishable (see also the right side of Figure \ref{fig-Example13}).

\begin{figure}[h!]
\centering
\def\svgwidth{300pt}
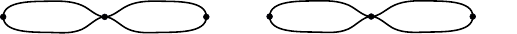
\caption{A non-labeled and a partially labeled graph.}
\label{fig-Example13}
\end{figure}
\end{Example}

\begin{Remark}\label{rem-ki1graphs} If $k_i=1$ for all $i$, then the valency condition implies that the genus at each vertex is $0$ and the vertices are $3$-valent. When forming the integral, the in the $z_i$ constant coefficient is just $1$, so we can neglect the $z_i$ and obtain Feynman integrals without vertex contributions in this case.
\end{Remark}

By Remarks \ref{rem-ki1covers} and \ref{rem-ki1graphs}, the equality in Theorem \ref{thm-mirror} specializes to the well-known relation involving the generating series of Hurwitz numbers and Feynman graphs, see, e.g., \cite[Theorem~9]{Dij95} and \cite[Theorem~2.6]{BBBM13}.

Using the correspondence Theorem \ref{thm-corres}, we can formulate a version of the mirror symmetry relation in Theorem \ref{thm-mirror}, where instead of the generating function of des-cendant Gromov--Witten invariants we use the generating function of tropical des-cendant Gromov--Witten invariants. It~turns out however that a finer version of a mirror symmetry relation naturally holds in the tropical world, which uses labeled tropical covers, multidegrees and refined Feynman integrals:

\begin{Definition}[labeled tropical cover]\label{def-LabeledTropCover}
Let $\pi$ be a tropical cover satisfying given psi-conditions with powers $k_1,\dots,k_n$ and denote the genus of a vertex of the source curve which is adjacent to end $i$ by $g_i$, where $g_i$ is given by $k_i$ via the psi-conditions (see Definition \ref{def-PsiAndPointConditions}). We can shrink the ends of the source curve and label the vertex that used to be adjacent to end $i$ with $x_i$. The cover $\pi$ is called \textit{labeled tropical cover} if there is an isomorphism of multigraphs sending a~Feynman graph $(\Gamma,\underline{g}')$ with a genus function (see Definition \ref{def-Feynmangraphs}) to the combinatorial type of the source curve of $\pi$, where the ends of the source curve are shrunk, such that $g_i'=g_i$ for all vertices. We say that $\pi$ is of type $\Gamma$.
\end{Definition}

Shortly, a labeled tropical cover is a tropical cover for which we label the vertices and edges of the source (vertices of different genus are distinguishable) according to a Feynman graph.

We fix a point $p_0\in E_{\TT}$. For a labeled tropical cover of $E_{\TT}$ of type $\Gamma$, we introduce its \emph{multidegree} as the vector $\underline{a}$ in $\NN^r$ with $k$-th entry $a_k=\big|\pi^{-1}(p_0)\cap q_k\big|\cdot \omega(q_k)$, where $\omega(q_k)$ denotes the expansion factor of the edge~$q_k$.
We define a labeled tropical descendant invariant
\begin{gather*}
\langle \tau_{k_1}(pt) \cdots \tau_{k_n}(pt) \rangle_{\Gamma,n}^{E,\underline{a},\trop}
\end{gather*}
as a count of labeled tropical covers of type $\Gamma$ and with multidegree $\underline{a}$ satisfying the prescribed point- and psi-conditions, again counted with multiplicity as in equation~ (\ref{eq-mult}). (Note that there are no nontrivial automorphism for a labeled tropical cover since all edges and vertices are distinguishable by their labeling.)

\begin{Theorem}[tropical mirror symmetry for $E_{\TT}$]\label{thm-refinedmirrortrop}
Fix $g\geq 2$, $n\geq 1$ and $k_1, \dots, k_n\geq 1$ satisfying $k_1+\dots+k_n = 2g-2$. Fix a Feynman graph $\Gamma$ such that the vertex $x_i$ has valency $k_i+2-2g_i$, and record the numbers $g_i$ in a genus vector $\underline{g}$.

Then we can express the series of descendant Gromov--Witten invariants of $E_{\TT}$ of type $\Gamma$ in terms of a Feynman integral:
\begin{gather*}
\sum_{\underline{a}\in\NN^r}\langle \tau_{k_1}(pt) \cdots \tau_{k_n}(pt) \rangle_{\Gamma,n}^{E,\underline{a},\trop}q_1^{a_1}\cdots q_r^{a_r} = I_{\Gamma,\underline{g}}(q_1,\dots,q_r).
\end{gather*}
\end{Theorem}
Theorem \ref{thm-mirror} follows from Theorem \ref{thm-refinedmirrortrop} using the correspondence Theorem \ref{thm-corres},
summing over all Feynman graphs $\Gamma$ such that $h^1(\Gamma)+\sum_{i=1}^n g_i=g$ (where the $g_i$ are defined by the valencies of the vertices as above), setting the $q_k$ equal to $q$ again for all $k$ and keeping track of automorphisms (as in the proof of Theorem 2.14 using Theorem 2.20 in \cite{BBBM13}).

We prove Theorem \ref{thm-refinedmirrortrop} in Section \ref{sec-bij} using Theorem \ref{thm-bij}, which establishes a bijection between labeled tropical covers contributing to a descendant Gromov--Witten invariant and monomials contributing to a term of the series used for the Feynman integrals.

\begin{Example}\label{ex-CalculationOfCoefficients}
Fix $g=2$. We want to use Theorem \ref{thm-refinedmirrortrop} to calculate contributions to $\langle \tau_{2}(pt) \tau_{0}(pt) \tau_{0}(pt) \rangle_{\Gamma,3}^{E,\underline{a},\trop}$ for two cases, where in the first case the covers contributing have a source curve with a nonzero genus function and in the second case the source curves have a~loop.

\medskip\noindent
\emph{First case:} we choose $a=(0,0,3)$ and $\Gamma$ as the left Feynman graph of Example~\ref{ex-underlyingFeynmanGraph}. So Theo\-rem~\ref{thm-refinedmirrortrop} tells us that we need to calculate the $q_1^0 q_2^0 q_3^3$-coefficient of $I_{\Gamma,\underline{g}}(q_1,q_2,q_3)$ with $\underline{g}=(1,0,0)$. We fix an order $\Omega$, namely the identity as we did in Example \ref{ex-stationarydescendant}. That is, we require that end $i$ is mapped to the point $p_i$.
Notice that the covers contributing to $\langle \tau_{2}(pt) \tau_{0}(pt) \tau_{0}(pt) \rangle_{\Gamma,3}^{E,(3,0,0),\trop}$ for $\Omega$ as above are the ones corresponding to the entries $(2,2)$ and $(3,1)$ in the table given in Example \ref{ex-stationarydescendant}. So we expect
\begin{align}\label{eq-ex-calculations1}
\frac{1}{24}+\frac{17\cdot 27}{24}=\frac{115}{6}
\end{align}
as the contribution.
We start by calculating terms of the propagators that contribute to the $q_1^0 q_2^0 q_3^3$-coefficient (we first let $w=1$ for $a_3$) in the product of the propagators such that their product is constant in the $x_i$, i.e., $l_1=l_2=l_3=0$,
\begin{gather*}
\tilde{P}\bigg(\frac{x_{1}}{x_{3}},q_3\bigg)=
\frac{4 \sinh \big( \frac{z_1}{2} \big) \sinh \big( \frac{z_3}{2}
 \big) \big( {\frac{ x_1}{x_3}}+{\frac{x_3}{x_1}} \big) {q_3}^{3}}{z_1 z_3}+\cdots,
\\
\tilde{P}\bigg(\frac{x_{2}}{x_{3}},q_2\bigg)=
\frac{4 \sinh \big( \frac{z_2}{2} \big) \sinh \big( \frac{z_3}{2} \big) x_2}{z_2 z_3 x_3}+\cdots,
\\
\tilde{P}\bigg(\frac{x_{1}}{x_{2}},q_1\bigg)=
\frac{4 \sinh \big( \frac{z_1}{2} \big) \sinh \big( \frac{z_2}{2} \big) x_1}{z_1 z_2 x_2}+\cdots.
\end{gather*}
Therefore,
\begin{gather*}
\Coef_{[q_1^{0}q_2^{0}q_3^{3}]} \frac{\tilde{P}(\frac{x_{1}}{x_{3}},q_3) \tilde{P}(\frac{x_{2}}{x_{3}},q_2) \tilde{P}(\frac{x_{1}}{x_{2}},q_1)} {\mathcal{S}(z_3)\mathcal{S}(z_2)\mathcal{S}(z_1)}
\\ \qquad
{}=
\frac {8 \sinh \big( \frac{z_1}{2} \big) \sinh \big( \frac{z_2}{2} \big)
\sinh \big( \frac{z_3}{2} \big) }{z_1 z_2 z_3}
= \dots + \frac{1}{1920} {z_1}^{4}+\frac{1}{576} {z_1}^{2}{z_2}^{2}+\frac{1}{576} {z_1}^{2} {z_3}^{2}
\\ \qquad\hphantom{=}
{}+\frac{1}{1920} {z_2}^{4}+ \frac{1}{576} {z_2}^{2} {z_3}^{2}
+ \frac{1}{1920} {z_3}^{4}+\frac{1}{24} {z_1}^{2}+\frac{1}{24} {z_2}^{2}+\frac{1}{24} {z_3}^{2}+1
\end{gather*}
and, hence, the $z_1^2z_2^0z_3^0$-coefficient is $\frac{1}{24}$, which is precisely the first summand of \eqref{eq-ex-calculations1}. The second summand is obtained by letting $w=3$ for $a_3$ such that
\begin{gather*}
\tilde{P}\bigg(\frac{x_{1}}{x_{3}},q_3\bigg)=
\frac{4 \sinh \big( \frac{3z_1}{2} \big) \sinh \big( \frac{3z_3}{2}
 \big) \big( {\frac{x_1^3}{x_3^3}}+{\frac{x_3^3}{x_1^3}} \big) {q_3}^{3}}{3 z_1 z_3}+\cdots,
\\
\tilde{P}\bigg(\frac{x_{2}}{x_{3}},q_2\bigg)=
\frac{4 \sinh \big( \frac{3z_2}{2} \big) \sinh \big( \frac{3z_3}{2} \big) x_2^3}{3z_2 z_3 x_3^3}+\cdots,
\\
\tilde{P}\bigg(\frac{x_{1}}{x_{2}},q_1\bigg)=
\frac{4 \sinh \big( \frac{3z_1}{2} \big) \sinh \big( \frac{3z_2}{2} \big) x_1^3}{3z_1 z_2 x_2^3}+\cdots
\end{gather*}
and therefore
\begin{gather*}
\Coef_{[q_1^{0}q_2^{0}q_3^{3}]} \frac{\tilde{P}(\frac{x_{1}}{x_{3}},q_3) \tilde{P}(\frac{x_{2}}{x_{3}},q_2) \tilde{P}(\frac{x_{1}}{x_{2}},q_1)} {\mathcal{S}(z_3)\mathcal{S}(z_2)\mathcal{S}(z_1)}
\\ \qquad
{}=\frac {8 \big(\sinh \big( \frac{3z_1}{2} \big)\big)^2 \big(\sinh \big( \frac{3z_2}{2} \big)\big)^2 \big(\sinh \big( \frac{3z_3}{2} \big)\big)^2 }{27 \sinh \big( \frac{z_1}{2} \big) \sinh \big( \frac{z_2}{2} \big) \sinh \big( \frac{z_3}{2} \big) z_1 z_2 z_3}
=\cdots{\frac {3369}{640}} {z_1}^{4}+{\frac {867}{64}}{z_1}^{2}
{z_2}^{2}+{\frac {867}{64}}{z_1}^{2}{z_3}^{2}
\\ \qquad \hphantom{=}
{}+{\frac {3369}{640}}{z_2}^{4}+{\frac {867}{64}}{z_2}^{2}{z_3}^{2}
+{\frac {3369}{640}}{z_3}^{4}+{\frac {153}{8}}{z_1}
^{2}+{\frac {153}{8}}{z_2}^{2}+{\frac {153}{8}}{z_3}^{2}+27,
\end{gather*}
where the $z_1^2z_2^0z_3^0$-coefficient is $\frac{153}{8}$ which equals the second summand of \eqref{eq-ex-calculations1}.

\medskip\noindent
\emph{Second case:} we choose $a=(2,0,0,1)$ and $\Gamma$ as the right Feynman graph of Example \ref{ex-underlyingFeynmanGraph}. By Theorem \ref{thm-refinedmirrortrop}, we need to calculate the $q_1^2 q_2^0 q_3^0 q_4^1$-coefficient of $I_{\Gamma,\underline{g}}(q_1,q_2,q_3,q_4)$ with $\underline{g}=0$. Again, we pick $\Omega$ as the order given by the identity.
As before, we calculate the terms of the propagators that contribute to the $q_1^2 q_2^0 q_3^0 q_4^1$-coefficient in the product of the propagators such that their product is constant in the $x_i$, i.e., $l_1=l_2=l_3=l_4=0$, and let $w=2$ for $a_1$, then
\begin{gather*}
\tilde{P}^{\lo}(q_1)=
\frac{2 (\sinh (z_1))^{2}{q_1}^{2}}{{z_1}^{2}} ,
\\
\tilde{P}\bigg(\frac{x_{1}}{x_{2}},q_2\bigg)=
{\frac {4\sinh \big( \frac{z_1}{2}\big) \sinh \big( \frac{z_2}{2} \big) x_1}{z_1 z_2 x_2}}+\cdots ,
\\
\tilde{P}\bigg(\frac{x_{2}}{x_{3}},q_3\bigg)=
{\frac {4\sinh \big( \frac{z_2}{2}\big) \sinh \big( \frac{z_3}{2} \big) x_2}{z_2 z_3 x_3}}+\cdots ,
\\
\tilde{P}\bigg(\frac{x_{1}}{x_{3}},q_4\bigg)=
\frac{4 \sinh \big( \frac{z_1}{2} \big) \sinh \big( \frac{z_3}{2}
 \big) \big( {\frac{x_1}{x_3}}+{\frac{x_3}{x_1}} \big) q_4}
{z_1 z_3}+\cdots
\end{gather*}
and
\begin{gather*}
\Coef_{[q_1^{2}q_2^{0}q_3^{0}q_4^{1}]} \frac{ \tilde{P}^{\lo}(q_1) \tilde{P}\big(\frac{x_{1}}{x_{3}},q_4\big) \tilde{P}\big(\frac{x_{2}}{x_{3}},q_3\big) \tilde{P}\big(\frac{x_{1}}{x_{2}},q_2\big)} {\mathcal{S}(z_3)\mathcal{S}(z_2)\mathcal{S}(z_1)}
\\ \qquad
{}=\frac { 16 \big( \sinh \big( z_1 \big) \big) ^{2}\sinh
 \big( \frac{z_1}{2} \big) \sinh \big( \frac{z_2}{2} \big) \sinh
 \big( \frac{z_3}{2} \big) }{{z_1}^{3} z_2 z_3}
=
2+\frac{3}{4}{z_1}^{2}+\frac{1}{12}{z_2}^{2}+\frac{1}{12}{z_3}^{2}+{\frac {113}{960}}{z_1}^{4}
\\ \qquad\hphantom{=}
{}+\frac{1}{32}{z_1}^{2}{z_2}^{2}
+ \frac{1}{32}{z_3}^{2} {z_1}^{2}+{\frac {1}{960}} {z_2}^{4}+{\frac {1}{288}}{z_3}^{2}{z_2}^{2}+{\frac {1}{960}}{z_3}^{4}+\cdots,
\end{gather*}
where the constant coefficient in the $z_i$ is $2$. If we let $w=1$ for $a_1$, we get $1$. This makes $3$ in total, which is the number we expect when using the table from Example~\ref{ex-stationarydescendant} again (entries $(1,1)$ and $(2,1)$).
\end{Example}

\subsection{The bijection}\label{sec-bij}
This subsection is devoted to the proof of the tropical mirror symmetry Theorem \ref{thm-refinedmirrortrop}. The main ingredient is a bijection of graph covers and monomials contributing to a Feynman integral.

Let $\Gamma$ be a Feynman graph (see Definition \ref{def-Feynmangraphs}). Fix a multidegree $\underline{a}$ and an order $\Omega$.

We can view $\Omega$ as an element in the symmetric group on $n$ elements, associating to $i$ the place $\Omega(i)$ that the vertex $x_i$ takes in the order $\Omega$.

Fix an orientation of $E_{\TT}$ and points $p_0,p_1,\dots,p_n$ ordered in this way when going around $E_{\TT}$ in the fixed orientation starting at $p_0$.

\begin{Definition}[graph covers of fixed order and multidegree] %\label{def-graphcover}
A \emph{graph cover} of type $\Gamma$, order~$\Omega$ and multidegree $\underline{a}$ is a (possibly leaky w.r.t.\ $(l_1,\dots,l_n)$, see Remark \ref{rem-leaking}) tropical cover $\pi$: $\Gamma'\rightarrow E_{\TT}$, where $\Gamma'$ is a metrization of $\Gamma$, such that the multidegree of $\pi$ at $p_0$ is $\underline{a}$ and such that $\pi^{-1}(p_{\Omega(i)})$ contains $x_i$. (Since there are $n$ point conditions and $n$ vertices, it follows that there is precisely one vertex of $\Gamma$ in each preimage $\pi^{-1}(p_{j})$).)

 We define $N^{l_1,\dots,l_n}_{\Gamma,\underline{a},\Omega}$ to be the weighted count of $(l_1,\dots,l_n)$-leaky graph covers of type $\Gamma$, order~$\Omega$ and multidegree $\underline{a}$, where we count each with the product of the expansion factors of the edges.

 If there is no mentioning of $l_1,\dots,l_n$, we refer to the case of no leaking as usual.
 \end{Definition}

Fix $g\geq 2$, $n\geq 1$ and $k_1, \dots, k_n\geq 1$ satisfying $k_1+\dots+k_n = 2g-2$.
 Let $\Gamma$ be a Feynman graph. Fix a multidegree $\underline{a}$ and an order $\Omega$. Assume that for each vertex $x_i$ of $\Gamma$, $k_i+\val(x_i)$ is even.

 \begin{Lemma}[graph covers and labeled tropical covers]\label{lem-graphtropcover}
There is a bijection between graph covers of type $\Gamma$, order $\Omega$ and multidegree $\underline{a}$ and labeled tropical covers $\pi\colon \Gamma'\rightarrow E_{\TT}$ contribu\-ting~to
\begin{gather*}
\langle \tau_{k_1}(pt) \cdots \tau_{k_n}(pt) \rangle_{\Gamma,n}^{E,\underline{a},\trop}
\end{gather*}
and satisfying $\pi(i)=p_{\Omega(i)}$.
 \end{Lemma}

\begin{proof}
Let $\pi\colon \Gamma'\rightarrow E_{\TT}$ be such a labeled tropical cover. We can describe the bijection as the map sending $\pi$ to a graph cover $\tilde{\pi}$ by shrinking marked ends of $\Gamma'$, labeling the vertex that used to be adjacent to end $i$ by $x_i$, and forgetting the genus at vertices. By definition of $\langle \tau_{k_1}(pt) \cdots \tau_{k_n}(pt) \rangle_{\Gamma,n}^{E,\underline{a},\trop}$, the graph cover is of type $\Gamma$. The multidegree is the same for the tropical cover and the graph cover. The set $\pi^{-1}(p_{\Omega(i)})$ contains $x_i$, since the marked end $i$ is mapped to $p_{\Omega(i)}$ by $\pi$. The inverse map associates the genus $\frac{k_i+2-\val(x_i)}{2}$ to the vertex $x_i$ (which is an integer by our assumption), and attaches the end marked $i$. Then the valence is $k_i+3-2g_i$ and the psi-condition is satisfied.
\end{proof}

\begin{Theorem}[bijection of graph covers and tuples in Feynman integrals]\label{thm-bij}
Let $\Gamma$ be a Feynman graph as in Definition $\ref{def-Feynmangraphs}$. Fix a multidegree $\underline{a}$ satisfying $a_k>0$ for all $k\leq s$, an order $\Omega$, and integers $l_1,\dots,l_n$.
As in Definition $\ref{def-Feynman}$, we use the notation $x_{k^1}$ and $x_{k^2}$ for the two vertices adjacent to the edge $q_k$, where we assume $x_{k^1}<x_{k^2}$ in $\Omega$.

There is a bijection between $(l_1,\dots,l_n)$-leaky graph covers of type $\Gamma$, order $\Omega$ and multidegree~$\underline{a}$, and tuples
\begin{equation}
\bigg((w_k)_{k=1,\dots,s},\bigg(\bigg(a_k, w_k\cdot\bigg(\frac{x_{k^i}}{x_{k^j}}\bigg)^{w_k}\bigg)\bigg)_{k=s+1,\dots,r}\bigg),
\label{eq-tuple}
\end{equation} where $i=1$
and $j=2$ if $a_k=0$, and $\{i,j\}=\{1,2\}$ otherwise, where $w_k$ divides $a_k$ if $a_k\neq 0$, and where the product of the fractions has exponent $l_i$ in $x_i$.

Moreover, the weighted count of graph covers equals the $q_1^{a_1}\cdots q_n^{a_n}$-coefficient of the refined Feynman integral:
\begin{gather*}
N^{l_1,\dots,l_n}_{\Gamma,\underline{a},\Omega}= \Coef_{[q_1^{a_1}\cdots q_r^{a_r}]} I^{l_1,\dots,l_n}_{\Gamma,\Omega}(q_1,\dots,q_r).
\end{gather*}
\end{Theorem}

\begin{Remark}[tuples and Feynman integrals]\label{rem-tuples}
Note that the products of second entries for $k>s$ of a tuple as in (\ref{eq-tuple}) are precisely the contributions showing up in the series
\begin{gather*}
\prod_{k=s+1}^{r} P\bigg(\frac{x_{k^1}}{x_{k^2}},q_k\bigg)
\end{gather*}
with the exponents of the $x_i$ given by the $l_i$, and the exponents of the $q_i$ given by the $a_i$. By~defi\-ni\-tion of the refined Feynman integral (see Definition \ref{def-Feynman}), adding a choice of summand~$w_k$ for each loop-edge $q_k$, $k\leq s$, each tuple contributes exactly $w_1\cdots w_r$ to the $q_1^{a_1}\cdots q_n^{a_n}$-coefficient of the refined Feynman integral $I^{l_1,\dots,l_n}_{\Gamma,\Omega}(q_1,\dots,q_r)$.

In particular, if $a_k=0$ for some $k\leq s$, the $q_1^{a_1}\cdots q_n^{a_n}$-coefficient of the refined Feynman integral $I^{l_1,\dots,l_n}_{\Gamma,\Omega}(q_1,\dots,q_r)$ is zero, and there are no tuples.
\end{Remark}

\begin{proof}[Proof of Theorem \ref{thm-bij}]
Given a tuple %$\Bigg(\Big(w_k\Big)_{k=1,\dots,s},\Big( \big(a_k, w_k\cdot\big(\frac{x_{k_i}}{x_{k_j}}\big)^{w_k}\big)\Big)_{k=s+1,\dots,r}\Bigg)$
as in (\ref{eq-tuple}), we construct a graph cover as follows. We~keep track of the cover by drawing the vertices and edges projecting onto their images. To~ease the drawing, we think of $E_{\TT}$ as being cut off at $p_0$ (see Example \ref{ex-relativestationarydescendant}).

We start by drawing vertices $x_i$ above the points $p_{\Omega(i)}$ in $E_{\TT}$.

For $k>s$ and for an entry $w_k\cdot\big(\frac{x_{k^i}}{x_{k^j}}\big)^{w_k}$, we draw an edge with expansion factor $w_k$ connecting vertex $x_{k^i}$ to vertex $x_{k^j}$. If $a_k=0$, we draw this edge in our cut picture direct, without passing over $p_0$. If $a_k>0$, we ``curl it'', passing over $p_0$ $\frac{a_k}{w_k}$ times before it reaches its end vertex.

We assume in our tuple that $i=1$ and $j=2$ if $a_k=0$. By Definition \ref{def-Feynman}, $x_{k^1}<x_{k^2}$ in $\Omega$, which implies that in our picture, the vertex $x_{k^1}$ is drawn before $x_{k^2}$ (in the orientation of $E_{\TT}$), which makes it possible to draw the edge $q_k$ directly without passing $p_0$.

Since $w_k$ divides $a_k$, it is possible to ``curl'' the edges $q_k$ with $a_k>0$ as required.

For $k\leq s$ and an entry $w_k$, we draw a loop-edge of weight $w_k$ connecting the vertex of $q_k$ to itself, ``curled'' over $p_0$ $\frac{a_k}{w_k}$ times.

In the drawing we created for the tuple (\ref{eq-tuple}), we have obviously drawn a graph cover with source curve of type $\Gamma$, since we connected the vertices $x_{k^1}$ and $x_{k^2}$ with the edge $q_k$. Furthermore, the multidegree is $\underline{a}$ because of our curling requirement. The order $\Omega$ is respected by the way we have drawn the vertices. It remains to show that the graph cover is $(l_1,\dots,l_n)$-leaky. To~see this, notice that the edges adjacent to vertex $x_i$ correspond to tuples whose fraction contains a power of $x_i$, and that the exponent equals $\pm$ the expansion factor of the edge, where the sign is positive if the edge leaves $x_i$ and negative if it enters $x_i$ (w.r.t.\ the orientation of $E_{\TT}$). Since we require the total power in $x_i$ to be $l_i$, the cover leaks $l_i$ at vertex $x_i$.

Clearly, the process can be reversed associating a tuple to a graph cover, and using the same arguments as before, the tuple satisfies the requirements from above. In particular, the entry $a_k$ of the multidegree of a cover with a loop-edge $q_k$ is nonzero. Thus, we have a bijection between graph covers and tuples.

The equality follows from Remark \ref{rem-tuples}, taking into account that a graph cover is counted with multiplicity the product of its expansion factors (which are the $w_i$) in $N^{l_1,\dots,l_n}_{\Gamma,\underline{a},\Omega}$.
\end{proof}

\begin{Example}
We illustrate the proof of Theorem \ref{thm-bij} by constructing a graph cover from a~given tuple as in \eqref{eq-tuple}. We let $\Gamma$ be the right graph of Example \ref{ex-underlyingFeynmanGraph}, $\Omega$ the identity, and $l_i=0$ for all $i$.
We choose the tuple
\begin{gather*}
\bigg( 1,\bigg( 2,1\cdot\bigg( \frac{x_1}{x_3} \bigg)^{-1} \bigg),\bigg( 0,1\cdot\bigg( \frac{x_1}{x_2} \bigg)^1 \bigg), \bigg( 0,1\cdot\bigg( \frac{x_2}{x_3} \bigg)^1 \bigg) \bigg)
\end{gather*}
and the order $\Omega$ given by the identity. Note that this tuple is not leaky. See Figure \ref{fig-Example3} for the following: We start by drawing the vertices $x_1$, $x_2$, $x_3$ above $p_1$, $p_2$, $p_3$. After that we add the non-curled edges $q_3$, $q_4$ which are given by the third and fourth entry of our tuple above. The edge $q_2$ is obtained by starting at $x_1$ and going left (we have a negative exponent in the second entry of our tuple), curling once (we want to pass $p_0$ twice with $q_2$) and ending at $x_3$. There is also one loop edge (the first entry of the tuple) adjacent to $x_1$ which does not curl. Since all weights of edges are $1$, we can also see from the graph that it is not leaky as we expected. The upper graph in Figure \ref{fig-Example3} inherits a metrization from downstairs. Thus we obtain a graph cover.

\begin{figure}[h!]\centering
\def\svgwidth{300pt}
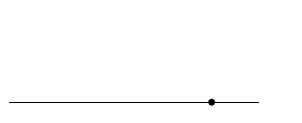
\caption{A graph cover constructed from a tuple. Note that this graph cover arises from cutting (and~labeling) the upper right source curve in Figure \ref{fig-Example1} along the preimages of a point $p_0$ (see also Example~\ref{ex-relativestationarydescendant}).}\label{fig-Example3}
\end{figure}
\end{Example}

As we have seen in Lemma \ref{lem-graphtropcover}, graph covers are closely related to tropical covers showing up in a tropical descendant Gromov--Witten invariant. However, the multiplicity of a tropical cover contains, besides the expansion factors for edges which already appear in the bijection in Theorem~\ref{thm-bij}, also factors for each vertex which can be viewed as local descendant Gromov--Witten invariants (see equation~(\ref{eq-localmult})).
The need to take these into account takes us to Feynman integrals with vertex contributions (see Definition~\ref{def-Feynmanvertex}).

We use Okounkov--Pandharipande's one-point series, i.e., the following nice form for the generating series of relative one-point descendant Gromov--Witten invariants, resp.\ for the local vertex multiplicities of tropical covers satisfying point and psi-conditions:
\begin{equation}
 \sum_{g\geq 0} \langle \mu| \tau_{2g-2+\ell(\mu)+\ell(\nu)}(pt)|\nu \rangle_{g, 1}^{\PP^1,d} \cdot z^{2g} = \frac{\prod\mathcal{S}(\mu_iz)\cdot \prod \mathcal{S}(\nu_iz)}{\mathcal{S}(z)}.\label{eq-onepointseries}
\end{equation}
Here, the $\mathcal{S}$-function is as in Definition \ref{def-prop}, and the product goes over all entries $\mu_i$ resp.\ $\nu_i$ of the two fixed partitions (see~\cite[Theorem 2]{OP06}, note that we consider Gromov--Witten invariants with the preimages of $0$ and $\infty$ marked).
Note that the $\mathcal{S}$-function is an even function (i.e., $\mathcal{S}(-z)=\mathcal{S}(z)$), since $\sinh$ is an odd function and quotients of odd functions are even.

\begin{proof}[Proof of Theorem \ref{thm-refinedmirrortrop}]
We prove the equality by restricting to the $q_1^{a_1}\cdots q_r^{a_r}$-coefficient on each side. It follows from Lemma \ref{lem-graphtropcover} that we can expand the left side as a sum over orders~$\Omega$, which we can do by definition of Feynman integral also on the right.
We thus have to show that the weighted count of labeled tropical covers contributing to $\langle \tau_{k_1}(pt) \cdots \tau_{k_n}(pt) \rangle_{\Gamma,n}^{E,\underline{a},\trop}$ and satisfying $\pi(i)=p_{\Omega(i)}$ equals $\Coef_{[q_1^{a_1}\cdots q_r^{a_r}]}I_{\Gamma,\underline{g},\Omega}(q_1,\dots,q_r)$.

To see this, note that by Remark \ref{rem-tuples} we deal with tuples as in Theorem \ref{thm-bij} when computing $\Coef_{[q_1^{a_1}\cdots q_r^{a_r}]}I_{\Gamma,\Omega}(q_1,\dots,q_r)$, however since we compute a Feynman integral with vertex contributions now each second entry $w_k\cdot\big(\frac{x_{k_i}}{x_{k_j}}\big)^{w_k}$ showing up in a tuple meets %an $\mathcal{S}$-function operator $ \mathcal{S}(z_{k^1} x_{k^1}\partial x_{k^1})\mathcal{S}(z_{k^2} x_{k^2}\partial x_{k^2})$ first, which ends up being
$ \mathcal{S}(w_k z_{k^1})\mathcal{S}(w_k z_{k^2})$ first.
By Theorem \ref{thm-bij}, the tuples are in bijection with graph covers. For a fixed graph cover corresponding to a fixed tuple, the vertex contributions in the Feynman integral thus produce factors of
$ \mathcal{S}(w_k z_{k^1})\mathcal{S}(w_k z_{k^2})$ for an edge of expansion factor $w_k$ connecting the vertices~$x_{k^1}$ and~$x_{k^2}$. Collecting those factors, sorting by $z_i$, and adding in the factor $\frac{1}{\mathcal{S}(z_i)}$ we have in the definition of Feynman integral with vertex contributions (see Definition \ref{def-Feynmanvertex}), we obtain for each vertex $x_i$ a contribution of
$\frac{\prod\mathcal{S}(\mu_jz_i)\cdot \prod \mathcal{S}(\nu_jz_i)}{\mathcal{S}(z_i)}$.
Here, the notation is set up as follows: we collect the expansion factors of all incoming edges adjacent to $x_i$ in the partition $\mu$ and those of all outgoing edges in the partition $\nu$. Taking the $z_i^{2g_i}$-coefficient, we obtain a local vertex contribu\-tion~of
\begin{gather*}
\langle \mu| \tau_{k_i}(pt)|\nu \rangle_{g_i, 1}^{\PP^1,|\mu|}
\end{gather*}
by the one-point series from equation~(\ref{eq-onepointseries}). By equation (\ref{eq-localmult}), this is exactly the local vertex multiplicity we need to take into account for the labeled tropical cover.
\end{proof}

\begin{Remark}
Let us compare the Tropical mirror symmetry Theorem \ref{thm-refinedmirrortrop} for descendant invariants with the version for Hurwitz numbers \cite[Theorem 2.20]{BBBM13}.
As we saw in Remarks~\ref{rem-ki1covers} and~\ref{rem-ki1graphs}, in the version for Hurwitz numbers, we only have to take $3$-valent graphs such that all vertices have genus $0$ into account.
Adding in descendants requires us to generalize in two ways: we need to include graphs whose vertices have other valencies, and whose vertices have genus.
The main ingredient in our proof of tropical mirror symmetry is the bijection between graph covers and monomials contributing to a Feynman integral, see Theorem \ref{thm-bij}. Graphs with vertices of valence bigger $3$ fit into this context.
The genus at vertices requires us to use local vertex multiplicities for the tropical covers, which are hard to translate to the Feynman integral world. The fact that the one-point series (\ref{eq-onepointseries}) can be expressed in a way separating contributions for the edges adjacent to a vertex makes it possible to incorporate these multiplicities in
a Feynman integral with vertex contributions as in Definition \ref{def-Feynmanvertex}.
\end{Remark}

\subsection{Quasimodularity}\label{sec-quasimod}

In the case that all $k_i=1$, the Mirror symmetry Theorem \ref{thm-mirror} specializes to the well-known relation involving the generating series of Hurwitz numbers and Feynman integrals for $3$-valent graphs without vertex contributions (see Remark \ref{rem-ki1graphs}).
This special case of the mirror symmetry relation was used in \cite{Dij95, KZ95} to prove that the generating function of Hurwitz numbers for $g\geq 2$ is a quasimodular form of weight $6g-6$.
Quasimodularity of generating functions of covers is a phenomenon studied beyond the case considered here, other important cases are generating functions of pillowcase covers \cite{EO06} or generating functions of numbers of covers of an elliptic curve with fixed ramifications with respect to the parity of the pullback of the trivial theta characteristic \cite{EOP08}.

Quasimodularity behaviour is desirable because it controls the asymptotic of the generating function. A series in $q$ is quasimodular if and only if it is in the polynomial ring generated by the three Eisenstein series $E_2$, $E_4$ and $E_6$ \cite{KZ95}. The \emph{weight} of a quasimodular form refers to its degree when viewed as a polynomial in the Eisenstein series. A series is called a quasimodular form of weight $w$ if it is a homogeneous polynomial of degree $w$ in the Eisenstein series, and it is called a quasimodular form of mixed weight if it is a non-homogeneous polynomial in the Eisenstein series.

From the Tropical mirror symmetry Theorem \ref{thm-refinedmirrortrop}, the generating function of descendant Gromov--Witten invariants of an elliptic curve obtains a natural stratification as sum over Feynman graphs, and, even finer, as sum over orders $\Omega$ for each Feynman graph (see Corollary \ref{cor-partiallyLabeled}).
If we fix a Feynman graph $\Gamma$ and a suitable genus function $\underline{g}$~-- if $k_i=1$ for all $i$, this means we fix a $3$-valent graph with genus $0$ at each vertex~-- we can study quasimodularity of individual summands.
We can consider summands $I_{\Gamma,\underline{g}}(q)$, or we can even break the sum into finer contributions by considering $I_{\Gamma,\underline{g},\Omega}(q)$ for a fixed order $\Omega$.

For the case that $k_i=1$ for all $i$, this study was initiated in \cite[Theorem 3.2]{BBBM13}, where it is conjectured that $I_{\Gamma,\underline{g}}(q)$ is quasimodular of weight $6g-6$. In \cite{GM16}, Goujard and M\"oller provide tools to study quasimodularity of generating series depending on Feynman graphs, and they prove that if $k_i=1$ for all $i$, each summand $I_{\Gamma,\underline{g},\Omega}(q)$ is a quasimodular form of mixed weight, where the highest appearing weight is $6g-6$. They also compute examples where lower weights appear. Since the whole sum (over all graphs, and over all orders) is quasimodular of weight $6g-6$, the contributions of lower weights must cancel in the sum. It remains an open question whether they already cancel in a summand $I_{\Gamma,\underline{g}}(q)$, i.e., when we sum over all orders $\Omega$, but for a fixed graph $\Gamma$.

In this section, we deduce from \cite[Theorem 6.1]{GM16} that $I_{\Gamma,\underline{g},\Omega}(q)$ is a quasimodular form of mixed weight also in the case of arbitrary $k_i$.
Also in the case of arbitrary $k_i$, quasimodularity (of mixed weight) of the whole generating series (the sum over all graphs, and all orders) was studied before \cite{Li16}. First, we interpret $I_{\Gamma,\underline{g},\Omega}(q)$ as a generating function of tropical covers:

\begin{Corollary}\label{cor-partiallyLabeled}
Fix $g\geq 2$, $n\geq 1$ and $k_1, \dots, k_n\geq 1$ satisfying $k_1+\dots+k_n = 2g-2$. Fix a~Feynman graph $\Gamma$ such that the vertex $x_i$ has valency $k_i+2-2g_i$ with $g_i\in \mathbb{N}$, and record the numbers $g_i$ in the genus function $\underline{g}$. Fix an order $\Omega$.
For $d\in \NN$, let $\langle \tau_{k_1}(pt) \cdots \tau_{k_n}(pt) \rangle_{\Gamma,n,\Omega}^{E,d,\trop}$ denote the number of (unlabeled) tropical covers (counted with multiplicity) contributing to $\langle \tau_{k_1}(pt) \cdots \tau_{k_n}(pt) \rangle_{g,n}^{E,d,\trop}$, for which the source curve has combinatorial type $\Gamma$ after shrinking the ends and satisfying $\pi(i)=p_{\Omega(i)}$.

Then we can express the generating function of these invariants in terms of the Feynman integral\vspace{-1ex}
\begin{gather*}
\sum_{d\in\NN}\langle \tau_{k_1}(pt) \cdots \tau_{k_n}(pt) \rangle_{\Gamma,n,\Omega}^{E,d,\trop}q^d =
\frac{1}{|\Aut(\ft_{\operatorname{edge}}(\Gamma),{\underline{g}})|} I_{\Gamma,\underline{g},\Omega}(q),
\end{gather*}
where $\ft_{\operatorname{edge}}$ is the map that forgets the edge labels of a Feynman graph, and automorphisms respect the remaining vertex labels and the genus function $($see Example $\ref{ex-automorphisms})$.
\end{Corollary}

\begin{proof}
Consider Theorem \ref{thm-refinedmirrortrop} and let $q_1=\dots=q_r=q$. With a similar argument as we use to deduce Theorem \ref{thm-mirror} from Theorem \ref{thm-refinedmirrortrop}, we also obtain an automorphism factor here. Fixing the order leads to labels on the vertices of the source curves, i.e., we need to consider automorphisms which respect partially labeled graphs as in Example \ref{ex-automorphisms}.
\end{proof}

\begin{Example}We want to express $I_{\Gamma,\underline{g},\Omega}(q)$ as polynomial in the Eisenstein series, where $\Omega$ is the identity, $\underline{g}=(0,0,0)$ or $\underline{g}=(1,0,0)$ and $\Gamma$ is any Feynman graph as shown in~Example~\ref{ex-underlyingFeynmanGraph}. So this example is a continuation of Examples \ref{ex-stationarydescendant} and \ref{ex-CalculationOfCoefficients}.

First, let $\Gamma_1$ be the left Feynman graph of Example \ref{ex-underlyingFeynmanGraph} and let $\underline{g}_1=(1,0,0)$. We calculate that\vspace{-1ex}
\begin{align*}
I_{\Gamma_1,\underline{g}_1,\Omega}(q)={}&
\frac{1}{20736}E_6(q)-\frac{1}{13824}E_2(q)E_4(q)+\frac{1}{41472}E_2^3(q)+\frac{1}{20736}E_4^2(q)
\\
& -\frac{1}{10368}E_2(q)E_6(q)+\frac{1}{20736}E_2^2(q)E_4(q)
\\
={}&\frac{1}{24}q+\frac{5}{2}q^2+\frac{39}{2}q^3+\frac{278}{3}q^4+\frac{1025}{4}q^5+738q^6 +\frac{4165}{3}q^7+3080q^8+\cdots.
\end{align*}
Notice that $I_{\Gamma_1,\underline{g}_1,\Omega}$ is of of mixed weight since $E_6$ and $E_2E_6$ are of different weight. Recall that we calculated $\frac{115}{6}$ as contribution to the $q^3$-coefficient. The other covers contributing are shown in Figure \ref{fig-Example1} of Example \ref{ex-stationarydescendant} and are the ones corresponding to the following entries of the table of Example \ref{ex-stationarydescendant}: $(2,3)$, $(2,4)$, $(2,5)$, $(3,2)$, $(3,3)$, $(4,1)$, $(4,2)$, $(4,3)$. Each of these covers contributes with $\frac{1}{24}$ such that in total we expect (see Corollary \ref{cor-partiallyLabeled})
{\samepage\begin{align*}
\Coef_{[q^3]}I_{\Gamma_1,\underline{g}_1,\Omega}(q)&=\frac{115}{6}+\frac{8}{24}
=\frac{39}{2},
\end{align*}
which matches our calculation.}

Second, we choose $\Gamma_2$ to be the right Feynman graph of Example \ref{ex-underlyingFeynmanGraph} and let $\underline{g}_2=0$. We~cal\-culate that
\begin{align*}
I_{\Gamma_2,\underline{g}_2,\Omega}(q)={}&
-\frac{1}{20736}E_6(q)+\frac{1}{13824}E_2(q)E_4(q)-\frac{1}{41472}E_2^3(q)+\frac{1}{20736}E_2(q)E_6(q)
\\
&-\frac{1}{13824}E_2^2(q)E_4(q)+\frac{1}{41472}E_2^4(q)
\\
={}&q^2+15q^3+76q^4+275q^5+720q^6+1666q^7+3440q^8+6129q^9+\cdots.
\end{align*}
Notice that, again, $I_{\Gamma_2,\underline{g}_2,\Omega}$ is of mixed weight, but $I_{\Gamma_1,\underline{g}_1,\Omega}+I_{\Gamma_2,\underline{g}_2,\Omega}$ is homogeneous. As above, we can verify the $q^3$-coefficient using Example \ref{ex-stationarydescendant}.

Third, we choose $\Gamma_3$ to be the middle Feynman graph of Example \ref{ex-underlyingFeynmanGraph} and let $\underline{g}_3=0$. In~this case, we obtain the homogeneous expression
\begin{align*}
I_{\Gamma_3,\underline{g}_3,\Omega}(q)&=
\frac{1}{20736}E_4^2(q)-\frac{1}{10368}E_2^2(q)E_4(q)+\frac{1}{20736}E_2^4(q)
\\&=
4q^2+48q^3+240q^4+800q^5+2160q^6+4704q^7+9920q^8+17280q^9+\cdots
\\&=
4\cdot(q^2+12q^3+60q^4+200q^5+540q^6+1176q^7+2480q^8+4320q^9+\cdots),
\end{align*}
where the factor $4$ in the last row is due to the automorphisms of the underlying Feynman graph, see Corollary \ref{cor-partiallyLabeled}. Again, we can verify the $q^3$-coefficient using Example~\ref{ex-stationarydescendant}.
\end{Example}

\begin{Corollary}%\label{cor-quasimod}
Fix $g$, $n$ and $k_1,\dots,k_n$ with $k_1+\dots+k_n = 2g-2$.
Let $\Gamma$ be a Feynman graph with $r$ edges $($see Definition $\ref{def-Feynmangraphs})$ and $\underline{g}$ a genus function satisfying $h^1(\Gamma)+\sum_{i=1}^n g_i=g$. Fix an order $\Omega$.
Then the Feynman integral
$I_{\Gamma,\underline{g},\Omega}(q)$~-- i.e., the generating function counting tropical covers for the tropical descendant Gromov--Witten invariant $\langle \tau_{k_1}(pt) \cdots \tau_{k_n}(pt) \rangle_{\Gamma,n}^{E,d,\trop}$ of type $\Gamma$ and order $\Omega$, see Corollary $\ref{cor-partiallyLabeled}$~-- is a quasimodular form of mixed weight, with highest occuring weight $2\cdot \big(r+\sum_{i=1}^n g_i\big)$.
\end{Corollary}

\begin{proof} This follows from \cite[Theorem~6.1]{GM16}, since the local vertex contributions we have to take into account for a vertex $x_i$ is polynomial of even degree $2g_i$ in the expansion factors of the adjacent edges by \cite[Theorem~4.1]{GM16} (see \cite{OP06, SSZ12}).
\end{proof}
This statement is essentially a byproduct of \cite[Corollary 6.2]{GM16} which states that the generating series of tropical covers with fixed ramification profiles (see \cite[Definition 2.1.3]{CJMR17}) and with fixed underlying graph $\Gamma$ and order $\Omega$ is a quasimodular form of mixed weight. The proof in \cite{GM16} detours by deducing the quasimodularity of the function above from the quasimodularity of our~$I_{\Gamma,\underline{g},\Omega}(q)$ (without explicitely stating this). The descendant Gromov--Witten invariants we focus on here are called Hurwitz numbers with completed cycles in \cite{GM16}, which is explained by the Okounkov--Pandharipande GW/H correspondence in \cite{OP06}, see also \cite{SSZ12}.

\section[Tropical mirror symmetry and the boson--fermion correspondence]
{Tropical mirror symmetry and the boson--fermion \\correspondence}\label{chap-fock}

The purpose of this section is to reveal the close relation between the proof of Theorem \ref{thm-mirror} in mathematical physics, using Fock spaces, and our tropical approach. Since the tropical setting requires a labeling of the underlying Feynman graphs and the use of the variables $q_1,\dots,q_r$ to distinguish degree contributions from the different edges, we enrich the Fock space approach by incorporating adequate labelings. This enlarges the set of operators, but makes it easier to distinguish contributions for a fixed Feynman graph to a matrix element. In this way, we extend the Fock space approach so that it gives an alternative proof of the tropical mirror symmetry Theorem \ref{thm-refinedmirrortrop}, which holds on a finer level. Our main ingredient is Theorem \ref{thm-fock}, proving the equality of the number of labeled tropical covers with fixed underlying source graph, fixed multidegree and order and a sum of matrix elements in a bosonic Fock space.

 For the sake of explicitness, we limit our considerations to the case of Hurwitz numbers, i.e., $k_i=1$ for all $i$, and we do not have vertex contributions for Feynman integrals (see Remarks~\ref{rem-ki1covers} and~\ref{rem-ki1graphs}). In particular, all our graphs are $3$-valent, have no loops and genus $0$ at vertices.
 Higher descendants resp.\ vertex contributions can be incorporated into our discussion also, but would increase the amount of notation largely~-- we would have to consider more summands for a bosonic vertex operator, and the tropical local vertex multiplicities would have to show up as coefficients of the bosonic vertex operator (see \cite[Section~5]{CJMR16}).

As shown in Figure~\ref{fig-Overview}, tropical geometry hands us a short-cut in the Fock space setting: we can relate the generating series of Hurwitz numbers directly to operators on the bosonic Fock space and do not need to invoke the fermionic Fock space and the boson--fermion correspondence, which is often viewed as the essence of mirror symmetry for elliptic curves.

\subsection{Hurwitz numbers as matrix elements}
We begin by shortly reviewing the bosonic Fock space approach for generating series of Hurwitz numbers.

The bosonic Heisenberg algebra $\mathcal{H}$ is the Lie algebra with basis $\alpha_n$ for $n\in \mathbb{Z}$ such that for $n\neq 0$ the following commutator relations are satisfied:
\begin{align*}
[\alpha_n, \alpha_m]=\left(n\cdot \delta_{n,-m}\right)\alpha_0,
\end{align*}
where $\delta_{n,-m}$ is the Kronecker symbol and $[\alpha_n, \alpha_m]:=\alpha_n\alpha_m-\alpha_m\alpha_n$.
The bosonic Fock space $F$ is a representation of $\mathcal{H}$. It is generated by a single ``vacuum vector'' $v_\varnothing$. The positive generators annihilate $v_\varnothing\colon  \alpha_n\cdot v_\varnothing=0$ for $n>0$, $\alpha_0$ acts as the identity and the negative operators act freely. That is, $F$ has a basis $b_\mu$ indexed by partitions, where
\begin{align*}
b_\mu=\alpha_{-\mu_1}\cdots \alpha_{-\mu_m}\cdot v_{\varnothing}.
\end{align*}

We define an inner product on $F$ by declaring $\langle v_\varnothing | v_\varnothing \rangle=1$ and $\alpha_n$ to be the adjoint of $\alpha_{-n}$.

We write
$\langle v|A|w\rangle$ for $\langle v|Aw\rangle$, where $v,w\in F$ and the operator $A$ is a product of elements in~$\mathcal{H}$, and $\langle A \rangle$ for $\langle v_\varnothing |A|v_\varnothing \rangle$. The first is called a \emph{matrix element}, the second a \emph{vacuum expectation}.

\begin{Definition} The \emph{cut-join} operator is defined by
\begin{equation*}
%\label{caj}M %=\sum_{\substack{i,j,k\in \mathbb{Z}\setminus\{0\}\\i+j+k=0}}\frac{1}{6} :\alpha_i \alpha_j \alpha_k:\;\;
M=\frac{1}{2} \sum_{k>0} \sum_{\substack{0< i, j \\i+j=k}} \alpha_{-j}\alpha_{-i}\alpha_k+\alpha_{-k}\alpha_{i}\alpha_{j}.
\end{equation*}
\end{Definition}

The relative invariants of $\mathbb{P}^1$ can be interpreted as a matrix element involving $M$ (notice that the invariants in questions are equal to double Hurwitz numbers by Okounkov--Pandharipande's GW/H correspondence, \cite[Theorem 1]{OP06}):

\begin{Proposition}\label{prop-p1fock} A relative Gromov--Witten invariant of $\PP^1$, resp.\ a double Hurwitz number, equals a matrix element on the bosonic Fock space:
\begin{gather*}
\langle \mu| \tau_{1}(pt)^n|\nu \rangle_{g, n}^{\PP^1,d,\bullet} = \frac{n!}{\prod_i\mu_i\cdot \prod_j \nu_j}\langle b_\mu | M^n | b_\nu\rangle.
\end{gather*}
\end{Proposition}
This statement follows by combining Wick's theorem with the correspondence Theorem \ref{thm-corresrel}: Wick's theorem (\cite[Theorem 5.4.3]{CJMR16}, \cite[Proposition 5.2]{BG14b}, \cite{Wic50}) expresses a matrix element as a~weighted count of graphs that are obtained by completing local pictures. It turns out that the graphs in question are exactly the tropical covers we enumerate to obtain $\langle \mu| \tau_{1}(pt)^n|\nu \rangle_{g, n}^{\PP^1,d, \trop}$, where $n!$ arises from fixing a set of points to which labeled ends are mapped to (rather than prescribing a point a labeled end should map to, see Definition \ref{def-PsiAndPointConditions}).

Notice that we have to use the disconnected theory here ($\bullet$), since the matrix element encodes \emph{all} graphs completing the local pictures and cannot distinguish connected and disconnected graphs.

The local pictures are built as follows: we draw one vertex for each cut-join operator. For an~$\alpha_n$ with $n>0$, we draw an edge germ of weight $n$ pointing to the right. If $n<0$, we draw an edge germ of weight $n$ pointing to the left.
For the two Fock space elements $b_\mu$ and $b_\nu$, we draw germs of ends: of weights $\mu_i$ on the left pointing to the right, of weights $\nu_i$ on the right pointing to the left. Wick's theorem states that the matrix element $\langle b_\mu | M^n | b_\nu\rangle$ equals a sum of graphs completing all possible local pictures, where each graph contributes the product of the weights of all its edges (including the ends).
A completion of the local pictures can be interpreted as a~tropical cover of $\PP^1_{\TT}$ (with suitable metrization).

\begin{figure}[h!]
\centering
\def\svgwidth{170pt}
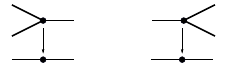
\caption{Local pictures of graphs with weights on the edges.}
\label{fig-Example4}
\end{figure}

The cut-join operator sums over all the possibilities of the local pictures for the graphs, i.e., it sums over all possibilities how a vertex of a tropical cover can look like (see Figure \ref{fig-Example4}).

\begin{Example}\label{ex-localpiecesWick}
Consider the local pieces shown below. There are three ways of completing them to a graph with local pictures like in Figure \ref{fig-Example4}.

\begin{figure}[h!]
\centering
\def\svgwidth{200pt}
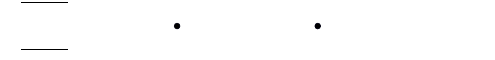
\end{figure}

\noindent The completed graphs are shown in Figure \ref{fig-Example6}. The product of the upper graph's edge weights (including the ends) is $12$, $4$ for the middle graph and $4$ for the lower graph. Hence Wick's theorem and Proposition \ref{prop-p1fock} yield
\begin{gather*}
\langle (2,1) | \tau_{1}(pt)^2| (2,1) \rangle_{2, 2}^{\PP^1,3,\bullet}=2! \cdot\left(3+1+\frac{1}{2}\right)=9,
\end{gather*}
where we have to divide the last summand by two because there is an automorphism exchanging the two edges that connect the vertices in the lower graph of Figure \ref{fig-Example6}.

\begin{figure}[h!]
\centering
\def\svgwidth{300pt}
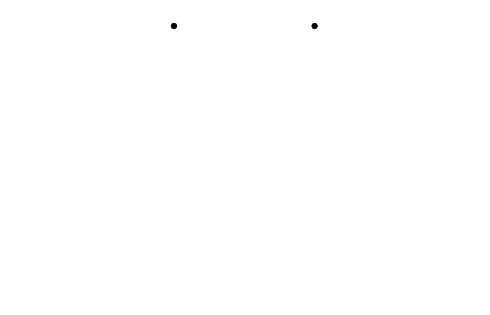
\caption{Completions of the local pieces above. Note that there is an automorphism exchanging the two bounded edges in the lower graph.}
\label{fig-Example6}
\end{figure}

\end{Example}

Combining Proposition \ref{prop-p1fock} with a degeneration argument, we can express Gromov--Witten invariants, resp.\ Hurwitz numbers of the elliptic curve in terms of matrix elements:

\begin{Proposition}\label{prop-etop1}
A Hurwitz number of the elliptic curve equals a weighted sum of double Hurwitz numbers of $\mathbb{P}^1$:
\begin{gather*}
\langle \tau_{1}(pt)^n \rangle_{g,n}^{E,d,\bullet}= \sum_{\mu\; \vdash d} \frac{\prod_i \mu_i}{|\Aut(\mu)|}
\langle \mu|\tau_{1}(pt)^n |\mu \rangle_{g-\ell(\mu), n}^{\PP^1,d,\bullet}.
\end{gather*}
Here, the sum goes over all partitions $\mu$ of $d$, $\mu_i$ denotes their entries, and $\ell(\mu)$ the length.
\end{Proposition}

\begin{Corollary}\label{cor-etop1}
A Hurwitz number of the elliptic curve $E$ equals a sum of matrix elements on the bosonic Fock space:
\begin{gather*}
\langle \tau_{1}(pt)^n \rangle_{g,n}^{E,d,\bullet} = \sum_{\mu\;\vdash d} \frac{n!}{|\Aut(\mu)|\prod_i\mu_i} \langle b_\mu | M^n | b_\mu\rangle.
\end{gather*}
\end{Corollary}

Proposition \ref{prop-etop1} is a corollary from the two correspondence Theorems \ref{thm-corres} and \ref{thm-corresrel}: given a~tropical cover of $E_{\TT}$, let $\mu$ be the partition encoding the weights of the edges mapping to the base point $p_0$. We mark the preimages of $p_0$, for which we have $|\Aut(\mu)|$ choices. For each choice, we cut off $E_{\TT}$ at $p_0$ and the covering curve at the preimages of $p_0$, obtaining a cover of $\PP^1_{\TT}$ with ramification profiles $\mu$ and $\mu$ above $\pm \infty$. The cut off tropical cover contributes to $\langle \mu|\tau_{1}(pt)^n |\mu \rangle_{g, n}^{\PP^1,d,\bullet, \trop}$, but its multiplicity differs from the multiplicity of the cover of $E_{\TT}$ by a~factor of $\prod \mu_i$, since the edges we cut off are no longer bounded.

\begin{Example}\label{ex-localpiecesWick2}
We want to calculate $\langle \tau_{1}(pt)^2 \rangle_{2,2}^{E,3,\bullet}$ using Corollary \ref{cor-etop1} and Wick's theorem. The partitions of $3$ are $(1,1,1)$, $(2,1)$ and $(3)$. The summand of $(2,1)$ follows from Example~\ref{ex-localpiecesWick}, namely $9\cdot 2=18$. The figure below shows how to complete the local pieces given by the partitions $(1,1,1)$ and $(3)$.

\begin{figure}[h!]
\centering
\def\svgwidth{230pt}
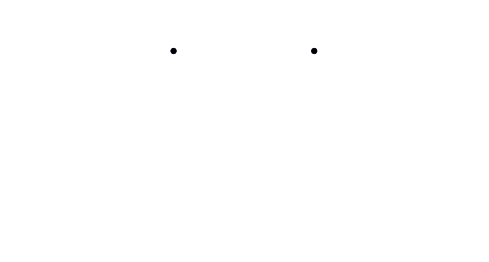
\caption{More completions of local pieces. Note that there are automorphisms of the upper graph that exchange the edges of weight one adjacent to a $3$-valent vertex.}
\label{fig-Example7}
\end{figure}

\noindent Note that there are in fact $9$ choices of how to complete the local pieces of $(1, 1, 1)$ since we can choose which ends (in the upper graph) the straight line should connect. Thus the upper graphs contribute ($9$ of them) $2!\cdot\frac{9\cdot 2}{4}\cdot \frac{1}{3!}=\frac{3}{2}$ and the lower graph contributes $2!\cdot 2\cdot 3=12$. Therefore $\langle \tau_{1}(pt)^2 \rangle_{2,2}^{E,3,\bullet}=\frac{63}{2}$.
\end{Example}

\subsection{Labeled matrix elements for labeled tropical covers}
Now we would like to link this Fock space language for Gromov--Witten invariants resp.\ Hurwitz numbers with tropical mirror symmetry. Recall that tropical mirror symmetry holds naturally on a fine level, giving an equality of the $q_1^{a_1}\cdots q_r^{a_r}$-coefficient of a Feynman integral $I_{\Gamma,\Omega}(q_1,\dots,q_r)$ and the number $N_{\Gamma,\underline{a},\Omega}$, which counts labeled covers of type $\Gamma$, with multidegree $\underline{a}$ and such that the order $\Omega$ is satisfied, i.e., the contribution to $\langle \tau_1(pt)^n\rangle_{\Gamma,n}^{E,\underline{a},\trop}$ of covers $\pi$ satisfying $\pi(i)=p_{\Omega(i)}$ (see Lemma \ref{lem-graphtropcover}, Theorems \ref{thm-bij} and~\ref{thm-refinedmirrortrop}).

Fix a Feynman graph $\Gamma$, a multidegree $\underline{a}$ and an order $\Omega$.
Remember that $\Gamma$ is a $3$-valent connected graph with first Betti number $g$, because of our restriction that $k_i=1$ for all $i$. In~particular, $\Gamma$ has no loops.

Our expression for $N_{\Gamma,\underline{a},\Omega}$ in terms of matrix elements (see Theorem \ref{thm-fock} below) involves a sum over all tuples $(w_k)_{k\colon a_k>0}$ with $w_k|a_k$ for all $k$ with $a_k>0$, since we incorporate the degeneration idea from above.

For a fixed choice of $(w_k)_k$, let $\Gamma'$ be the graph obtained from $\Gamma$ by cutting the edge $q_k$ exactly~$\frac{a_k}{w_k}$ times. We introduce the following labels for the (cut) edges of $\Gamma'$: we denote the pieces by $q_{k,1},\dots,q_{k,\frac{a_k}{w_k}+1}$. There are at most $a_k+1$ pieces, depending on $w_k$. For an edge which is not cut, i.e., $a_k=0$, we call it $q_{k,1}$ to consistently have two indices for the edge labels in $\Gamma'$.

We enlarge our set of operators in a way that allows to distinguish the edges of the cut graph $\Gamma'$:
Let the $\alpha^{k,j}_n$, for each $k=1,\dots,r$, $j=1,\dots,a_k+1$, and $n\in \mathbb{Z}\setminus\{0\}$, satisfy the commutator relations
\begin{align*}
 [\alpha^{k,j}_n, \alpha^{l,i}_m]=\left(n\cdot\delta_{k,l}\cdot \delta_{j,i} \cdot\delta_{n,-m}\right)\alpha_0.
\end{align*}
As before, we let the bosonic Fock space $F$ be generated by $v_\varnothing$, following the rules from before: $\alpha^{k,j}_n\cdot v_\varnothing=0$ for $n>0$, $\alpha_0$ acts as identity, and the operators with negative subscript act freely.
We let $\langle v_\varnothing | v_\varnothing \rangle=1$ and let $\alpha^{k,j}_n$ be the adjoint of $\alpha^{k,j}_{-n}$.

\begin{Definition} Let $\Gamma$, $\underline{a}$ and $(w_k)_k$ be as above. Let $x_i$ be a vertex of $\Gamma$. We denote the three adjacent edges by $q_{i_1}$, $q_{i_2}$ and $q_{i_3}$.
For $l=1,2,3$ set $c_l=\frac{a_{i_l}}{w_{i_l}}+1$ if $a_{i_l}>0$ and $c_l=1$ else. We~also set $d_{m_l}=c_l$ if $m_l>0$ and $d_{m_l}=1$ otherwise.

The labeled cut-join operator for the vertex $x_i$ is
\begin{align*}
M_{i}=&\sum_{\substack{m_1,m_2,m_3\in \mathbb{Z}\setminus\{0\}\\m_1+m_2+m_3=0}} \alpha^{i_1,d_{m_1}}_{m_1} \alpha^{i_2,d_{m_2}}_{m_2}\alpha^{i_3,d_{m_3}}_{m_3}.
\end{align*}
Since the first superscript differs for the $\alpha$-operators in a summand, the commutator relations imply that these factors can be permuted within a summand without changing the cut-join operator.
\end{Definition}

\begin{figure}[h!]
\centering
\def\svgwidth{270pt}
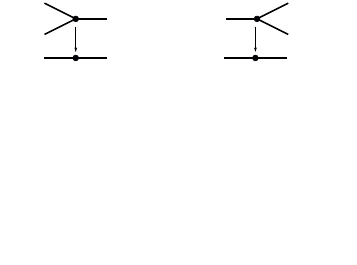
\caption{Local pictures of graphs with weights $m_1$, $m_2$, $m_3$ on the labeled edges $q_{i_1}$, $q_{i_2}$, $q_{i_3}$.}
\label{fig-Example8}
\end{figure}

This operator sums over all possibilities of how, locally, a vertex with its adjacent edge germs can be arranged, as shown in Figure \ref{fig-Example8}.

Notice that, compared to the (unlabeled) cut-join operator, we do not need a factor of $\frac{1}{2}$ which had to be there to take automorphisms into account resp.\ to undo overcounting by distinguishing edges which are not distinguishable. Here, all edges are labeled and thus distingushable.

\begin{Theorem}\label{thm-fock}
For a fixed Feynman graph $\Gamma$, multidegree $\underline{a}$ and order $\Omega$, the count of labeled tropical covers of $E_{\TT}$ of type $\Gamma$ and of the right multidegree and order $($see Lemma $\ref{lem-graphtropcover}$, Theorems~$\ref{thm-bij}$ and~$\ref{thm-refinedmirrortrop})$ equals a sum of matrix elements:
\begin{gather*}
N_{\Gamma,\underline{a},\Omega}=\sum_{\substack{(w_k)_k\\w_k|a_k}}
\prod_{k=1}^r \bigg(\frac{1}{w_k}\bigg)^{\frac{a_k}{w_k}}\cdot
\Biggl\langle \prod_{k=1}^r \prod_{l=1}^{\frac{a_k}{w_k}} \alpha^{k,l}_{-w_k} v_{\varnothing} \Biggl| \prod_{i=1}^n M_{\Omega^{-1}(i)} \Biggl| \prod_{k=1}^r \prod_{l=2}^{\frac{a_k}{w_k}+1} \alpha^{k,l}_{-w_k} v_{\varnothing} \Biggl\rangle.
\end{gather*}
\end{Theorem}

\begin{proof}
We use Wick's theorem: the right hand side is a sum over all possible ways to combine the local pictures given by the cut-join operators to a graph $\Gamma'$ that covers $\PP^1_{\TT}$. Our local pictures are now vertices with labels~$x_{\Omega^{-1}(i)}$. For each $\alpha_n^{k,j}$, we have an adjacent edge germ with label $q_{k,j}$ of weight $|n|$, pointing to the right if $n$ is positive and to the left otherwise.
Fix a graph~$\Gamma'$ which is a completion of such local pictures. The preimages of $-\infty$ are leaf vertices of~$\Gamma'$ whose adjacent edges are labeled $q_{k,1},\dots,q_{k,\frac{a_k}{w_k}}$ and are of weight $w_k$ (for all $k$). The preimages of $\infty$ are leaf vertices whose adjacent edges have labels $q_{k,2},\dots,q_{k,\frac{a_k}{w_k}+1}$, also of weight $w_k$. Since the $\alpha$-operators in the cut-join operator only have the values $1$ or $\frac{a_k}{w_k}+1$ as their second superscript, the commutator relation guarantees that the leaves of $q_{k,2},\dots,q_{k,\frac{a_k}{w_k}}$ over $-\infty$ have to be connected to the leaves with the corresponding label over $\infty$.
The leaf adjacent to $q_{k,1}$ over $-\infty$ is merged with an interior vertex adjacent to $q_k$, by definition of the labeled cut-join operator which depends on $\Gamma$. The same holds for the leaf adjacent to $q_{k,\frac{a_k}{w_k}+1}$.

To produce a tropical cover of $E_\TT$, we glue $\Gamma'$ as follows: for all $k$ and for $i=1,\dots,\frac{a_k}{w_k}$, the leaf of $q_{k,i}$ over $-\infty$ is attached to the leaf of $q_{k,i+1}$ over $\infty$. Identifying the edges $q_{k_1},\dots,q_{k,\frac{a_k}{w_k}+1}$ (which are subdivided by $2$-valent vertices obtained from gluing) to one edge $q_k$, we obtain a~graph cover of $E_{\TT}$ of type $\Gamma$ which is of the right order and multidegree: the order is imposed by the order in which we multiply the cut-join operators, the multidegree is given by the ``curls'' of the edge $q_k$, which has weight $w_k$ and which is curled $\frac{a_k}{w_k}$ times by our way of gluing.

Obviously, each tropical cover of type $\Gamma$ and multidegree $\underline{a}$ with order $\Omega$ can be obtained by gluing a graph $\Gamma'$ that arises with Wick's theorem from the right hand side.

 On the right hand side, a graph $\Gamma'$ that we produce with Wick's theorem contributes with the product of the weights of all edges which are connected, including the ends. For an edge $q_k$ (which is cut into $w_k^{\frac{a_k}{w_k}+1}$ pieces in $\Gamma'$) with $a_k>0$, we thus obtain a factor of $w_k^{\frac{a_k}{w_k}+1}$, where we actually only want $w_k$ for the tropical multiplicity. This is taken care of by the pre-factor before the summands on the right.
\end{proof}

\begin{Example}
Fix the multidegree $\underline{a}=(2,1,0)$, an order $\Omega$ on the vertices $x_1$, $x_2$ such that $x_1<x_2$ and the following Feynman graph:
\begin{figure}[h!]
\centering
\def\svgwidth{90pt}
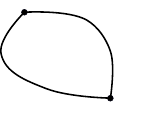
%\label{fig-Example9}
\end{figure}
\noindent Fix points $p_1$, $p_2$ on $E_{\TT}$. We obtain a labeled tropical cover of $\PP_{\TT}^1$ that can be glued to a cover of $E_{\TT}$ of type $\Gamma$ by choosing local pieces (see Figure \ref{fig-Example8}).

Notice that there are two choices of the expansion factor $w_1$, namely $w_1=1$ or $w_1=2$. We~start with $w_1=2$ and obtain the following source curve of a tropical cover to $\PP_{\TT}^1$, where the local pieces are indicated by boxes.
\begin{figure}[h!]
\centering
\def\svgwidth{380pt}
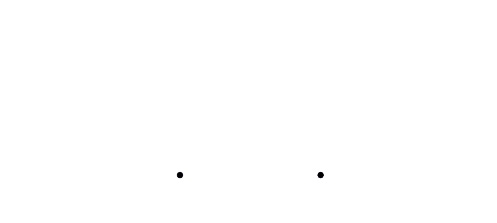
%\label{fig-Example10}
\end{figure}
\noindent In case of $w_1=2$, there are no other choices of local pieces that fit $\Gamma$ than the ones shown above. If we choose $w_1=1$, then another valid choice of local pieces is shown below.

\begin{figure}[h!]
\centering
\def\svgwidth{380pt}
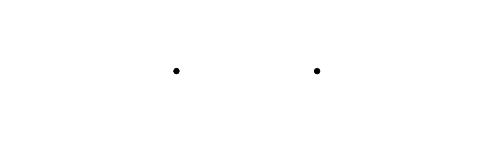
%\label{fig-Example11}
\end{figure}

Note that these two graphs are labeled version of the middle and upper graph of Figure \ref{fig-Example6} and the upper graph of Figure \ref{fig-Example7}. However, we do not get all graphs of Examples \ref{ex-localpiecesWick} and \ref{ex-localpiecesWick2} since we fixed $\Gamma$ and $\underline{a}$.

\end{Example}

\subsection{From matrix elements to Feynman integrals}
Finally, we link the matrix elements on the right of the equation in Theorem \ref{thm-fock} with Feynman integrals.

For this purpose, we introduce formal variables for our vertices in the labeled cut-join operators:

\begin{Definition} Let $\Gamma$, $\underline{a}$ and $(w_k)_k$ be as above. For $l=1,2,3$ set $c_l=\frac{a_{i_l}}{w_{i_l}}+1$ if $a_{i_l}>0$ and $c_l=1$ otherwise.
We also set $d_{m_l}=c_l$ if $m_l>0$ and $d_{m_l}=1$ otherwise.

The labeled cut-join operator for the vertex $x_i$ is
\begin{gather*}
M(x_i)= \sum_{m_1,m_2,m_3\in \mathbb{Z}\setminus\{0\}} \alpha^{i_1,d_{m_1}}_{m_1} x_i^{m_1} \alpha^{i_2, d_{m_2}}_{m_2}x_i^{m_2} \alpha^{i_3, d_{m_3}}_{m_3}x_i^{m_3}.
\end{gather*}
\end{Definition}
Here, we treat the cut-join operators as formal series in $x_1,\dots,x_n$.

With this, we can rewrite the equation of Theorem \ref{thm-fock} as follows:
\begin{gather}
N_{\Gamma,\underline{a},\Omega}=\Coef_{x_1^0\cdots x_n^0}\sum_{\substack{(w_k)_k\\w_k|a_k}}
\prod_{k=1}^r \bigg(\frac{1}{w_k}\bigg)^{\frac{a_k}{w_k}}\nonumber
\\ \hphantom{N_{\Gamma,\underline{a},\Omega}=}
{}\times\Biggl\langle \prod_{k=1}^r \prod_{l=1}^{\frac{a_k}{w_k}} \alpha^{k,l}_{-w_k} v_{\varnothing} \Biggl| \prod_{i=1}^n M(x_{\Omega^{-1}(i)}) \Biggl| \prod_{k=1}^r \prod_{l=2}^{\frac{a_k}{w_k}+1} \alpha^{k,l}_{-w_k} v_{\varnothing} \Biggl\rangle.
\label{eq-fockwithx}
\end{gather}

Each matrix element on the right hand side is now a series in $x_1,\dots,x_n$ when evaluated.

\begin{Lemma} \label{lem-fockfeynman}
Fix $\Gamma$, $\underline{a}$ and $\Omega$ as above.
The matrix elements of equation~$(\ref{eq-fockwithx})$, viewed as series in $x_1,\dots,x_n$ equals the following product:
\begin{gather*}
\sum_{\substack{(w_k)_k\\w_k|a_k}}
\prod_{k=1}^r \bigg(\frac{1}{w_k}\bigg)^{\frac{a_k}{w_k}}\cdot
\Biggl\langle \prod_{k=1}^r \prod_{l=1}^{\frac{a_k}{w_k}} \alpha^{k,l}_{-w_k} v_{\varnothing} \Biggl| \prod_{i=1}^n M(x_{\Omega^{-1}(i)}) \Biggl| \prod_{k=1}^r \prod_{l=2}^{\frac{a_k}{w_k}+1} \alpha^{k,l}_{-w_k} v_{\varnothing} \Biggl\rangle
\\ \qquad
{}= \prod_{k\colon a_k>0} w_k \cdot \Biggl(\bigg(\frac{x_{k^1}}{x_{k^2}}\bigg)^{w_k}+\bigg(\frac{x_{k^2}}{x_{k^1}}\bigg)^{w_k} \Biggl) \cdot \prod_{k\colon a_k=0} \Biggl( \sum_{w_k>0} w_k\cdot\bigg(\frac{x_{k^1}}{x_{k^2}}\bigg)^{w_k}
\Biggl).
\end{gather*}
Here, $x_{k^1}$ and $x_{k^2}$ denote the vertices adjacent to the edge $q_k$, where in the order $\Omega$ we have $x_{k^1}<x_{k^2}$.
\end{Lemma}

\begin{proof}
Let $q_k$ be an edge with $a_k=0$.
Since $a_k=0$, an $\alpha$ with first superscript $k$ does not show up in the vectors of the matrix element, only in the labeled cut-join operators. Also, the second superscript must be $1$, and it appears for exactly two cut-join operators, namely the one for $x_{k^1}$ and the one for $x_{k^2}$. Thus, we draw an edge germ labeled $q_{k,1}$ at $x_{k^1}$ and an edge germ labeled $q_{k,1}$ at $x_{k^2}$. These are the only edge germs with this label.
To obtain a nonzero contribution to the matrix element, the edge germ at $x_{k^1}$ must point to the right and the one at $x_{k^2}$ must point to the left. Furthermore, they must have the same weight $w_k$.
There is no restriction on the weight $w_k$. (The balancing condition is imposed only after we take the $x_1^0\cdots x_n^0$-coefficient in equation~(\ref{eq-fockwithx}).)
So, for any $w_k>0$, we have nonzero contributions to the matrix elements above with an $\alpha_{w_k}^{k,1}$ in the cut-join operator $M(x_{k^1})$ and an $\alpha_{-w_k}^{k,1}$ in the cut-join operator $M(x_{k^2})$. Combining those $\alpha$-operators with the respective power of the variable, we obtain $\alpha_{w_k}^{k,1} \cdot x_{k^1}^{w_k} \cdot \alpha_{-w_k}^{k,1}\cdot x_{k^2}^{-w_k}$, which, after applying the commutator relation and simplifying becomes $w_k\cdot \big(\frac{x_{k^1}}{x_{k^2}}\big)^{w_k}$.

We have treated the sum of matrix elements as a weighted sum of graphs. Any nonzero summand must have an edge connecting the edge germs above, and it can be of any weight. More precisely, if we have a graph with such an edge of a certain weight, we also have all summands that correspond to the same graph, but with the weight of the edge varying.
Thus, we can pull out a factor
\begin{gather*}
\sum_{w_k>0} w_k\cdot \bigg(\frac{x_{k^1}}{x_{k^2}}\bigg)^{w_k}
\end{gather*}
for the edge $q_k$.

Let us now consider an edge $q_k$ with $a_k>0$. The matrix elements on the left are summed over all $w_k|a_k$.
For the local pictures, we draw end germs of weight $w_k$ on the left pointing to the right, with labels $q_{k,1},\dots,q_{k,\frac{a_k}{w_k}}$, and on the right, pointing to the left, with labels $q_{k,2},\dots,q_{k,\frac{a_k}{w_k}+1}$. We~use the commutator relations for the $\alpha$ in charge of connecting the ``curls'' $q_{k,2},\dots,q_{k,\frac{a_k}{w_k}}$, they produce a factor of $w_k$ which is cancelled by the pre-factor. The edge germ $q_{k,1}$ must be connected to an edge germ appearing in a cut-join operator, that can be either $M(x_{k^1})$ or~$M(x_{k^2})$. The edge germ $q_{k,\frac{a_k}{w_k}+1}$ must also be connected to an edge germ of a cut-join operator, necessarily the other on in the choice of $M(x_{k^1})$ or $M(x_{k^2})$. Thus, we either have
\begin{gather*}
\alpha_{w_k}^{k,1}\cdot \alpha_{-w_k}^{k,1}\cdot x_{k^1}^{-w_k}\cdot \alpha_{w_k}^{k,\frac{a_k}{w_k}+1}\cdot x_{k^2}^{w_k}\cdot \alpha_{-w_k}^{k,\frac{a_k}{w_k}+1}
\end{gather*}
or
\begin{gather*}
\alpha_{w_k}^{k,1}\cdot \alpha_{-w_k}^{k,1}\cdot x_{k^2}^{-w_k}\cdot \alpha_{w_k}^{k,\frac{a_k}{w_k}+1}\cdot x_{k^1}^{w_k}\cdot \alpha_{w_k}^{k,\frac{a_k}{w_k}+1}
\end{gather*}
(notice the subscript changes sign when we let factors jump in the scalar product, by convention of the adjoints). Taking the commutator relations into account, and realizing that one factor of~$w_k$ is again cancelled by the pre-factor, we obtain either $w_k\cdot \big(\frac{x_{k^1}}{x_{k^2}}\big)^{w_k}$ or $w_k\cdot \big(\frac{x_{k^2}}{x_{k^1}}\big)^{w_k}$. For any graph which produces the first factor, we can connect the edge germs differently thus obtaining the graph which produces the second factor. Also, $w_k$ was imposed by the summand we picked on the left hand side. But for a given graph with an edge of weight $w_k$, we also have the analogous graph (with fewer or more curls), where the edge has another weight which divides~$a_k$. Thus, for the edge $q_k$ we obtain a total factor of
\begin{gather*}
\sum_{w_k|a_k} w_k\cdot \bigg(\bigg(\frac{x_{k^1}}{x_{k^2}}\bigg)^{w_k}+\bigg(\frac{x_{k^2}}{x_{k^1}}\bigg)^{w_k}\bigg).
\tag*{\qed}
\end{gather*}
\renewcommand{\qed}{}
\end{proof}

With this, we can now give an alternative proof of Theorem \ref{thm-refinedmirrortrop} (in the case $k_i=1$ for all $i$), which follows the more traditional Fock space approach, now with a larger set of operators in charge of the labels. In the tropical world, we can however take a shortcut avoiding the fermionic Fock space and relying on Wick's theorem instead.

\begin{proof}[Proof of Theorem~\ref{thm-refinedmirrortrop} in the case $\boldsymbol{k_i=1}$ for all $\boldsymbol i$]
We prove the equality by restricting to the $q_1^{a_1}\cdots q_r^{a_r}$-coefficient on each side. It follows from Lemma \ref{lem-graphtropcover} that we can expand the left side as a sum over orders $\Omega$, which we can do by definition of Feynman integral also on the right.
We thus have to show that the weighted count $N_{\Gamma,\underline{a},\Omega}$ of labeled tropical covers contributing to $\langle \tau_{k_1}(pt) \cdots \tau_{k_n}(pt) \rangle_{\Gamma,n}^{E,\underline{a},\trop}$ and satisfying $\pi(i)=p_{\Omega(i)}$ equals $\Coef_{[q_1^{a_1}\cdots q_r^{l_r}]}I_{\Gamma,\underline{g}\Omega}(q_1,\dots,q_r)$.
By~definition of a Feynman integral (see Definition \ref{def-Feynman}), the $x_1^0\cdots x_n^0$-coefficient of the right hand side of Lemma~\ref{lem-fockfeynman} equals the $q_1^{a_1}\cdots q_r^{a_r}$-coefficient of $I_{\Gamma,\Omega}(q_1,\dots,q_r)$. Using Lemma \ref{lem-fockfeynman} and equation~(\ref{eq-fockwithx}) (which follows from Theorem~\ref{thm-fock} relying on Wick's theorem), it follows that it also equals $N_{\Gamma,\underline{a},\Omega}$. The statement is proved.
\end{proof}

\subsection*{Acknowledgements}
We would like to thank Renzo Cavalieri, Elise Goujard, Gerhard Hiss, Martin M\"oller and Dhruv Ranganathan for helpful discussions.
Gefördert durch die Deutsche Forschungsgemeinschaft (DFG) -- Projektnummer 286237555 -- TRR~195 [Funded by the Deutsche Forschungsgemeinschaft (DFG, German Research Foundation) -- Project ID 286237555 -- TRR~195].
The authors have been supported by Project I.10 (INST 248/237-1) of TRR 195. Computations have been made with \textsc{Singular} using the \textit{ellipticcovers} library.
Part of this work was completed during the Mittag-Leffler programm \emph{Tropical geometry, amoebas and polytopes} in spring~2018. The~authors would like to thank the institute for hospitality and excellent working conditions.
We~would like to thank the anonymous referees for helpful suggestions to improve the paper.

\pdfbookmark[1]{References}{ref}
\LastPageEnding

\end{document}